\documentclass[a4paper,reqno, oneside]{amsart}
\usepackage[english]{babel}
\usepackage{amsmath, amssymb, amsthm, amscd}
\usepackage{enumerate}
\usepackage{palatino}
\usepackage{mathpazo}
\usepackage{paralist}
\usepackage[a4paper]{geometry}
\usepackage[all]{xy}
\usepackage{graphicx}

\usepackage{tikz} 
\usepackage{tikz-cd}

\usepackage{csquotes}

\usepackage{mathtools}

\DeclareMathOperator{\id}{id}

\newcommand{\NN}{\mathbb{N}}
\newcommand{\ZZ}{\mathbb{Z}}
\newcommand{\RR}{\mathbb{R}}
\newcommand{\CC}{\mathbb{C}}

\newcommand{\Ad}{\mathrm{Ad}}

\DeclareMathOperator{\Exp}{Exp}

\renewcommand{\setminus}{\smallsetminus}

\DeclareMathOperator{\cpr}{\vee}

\makeatletter
\newcommand{\ostar}{\mathbin{\mathpalette\make@circled\ast}}
\newcommand{\make@circled}[2]{%
  \ooalign{$\m@th#1\smallbigcirc{#1}$\cr\hidewidth$\m@th#1#2$\hidewidth\cr}%
}
\newcommand{\smallbigcirc}[1]{%
  \vcenter{\hbox{\scalebox{0.77778}{$\m@th#1\bigcirc$}}}%
}
\makeatother

\theoremstyle{plain}
\newtheorem{theorem}{Theorem}[section]
\newtheorem{prop}[theorem]{Proposition}
\newtheorem{lemma}[theorem]{Lemma}
\newtheorem{cor}[theorem]{Corollary}  

\newtheorem*{thm}{Theorem}

\theoremstyle{definition}
\newtheorem{definition}[theorem]{Definition}

\theoremstyle{remark}
\newtheorem{remark}[theorem]{Remark}

\numberwithin{equation}{section}

\begin{document}
\setlength{\parindent}{0.cm}

\title{On the string topology of symmetric spaces of higher rank}

\author{Philippe Kupper}
\address{Institut für Algebra und Geometrie, KIT Karlsruhe, 
Englerstrasse 2, 76131 Karlsruhe}
\email{philippe.kupper@kit.edu}

\author{Maximilian Stegemeyer}
\address{Mathematisches Institut, Universit\"at Freiburg, Ernst-Zermelo-Straße 1, 79104 Freiburg, Germany}
\email{maximilian.stegemeyer@math.uni-freiburg.de}
\date{\today}
\maketitle

\begin{abstract}
The homology of the free and the based loop space of a compact globally symmetric space can be studied through explicit cycles.
We use cycles constructed by Bott and Samelson and by Ziller to study the string topology coproduct and the Chas-Sullivan product on compact symmetric spaces.
We show that the Chas-Sullivan product for compact symmetric spaces is highly non-trivial for any rank and we prove that there are many non-nilpotent classes whose powers correspond to the iteration of closed geodesics.
Moreover, we show that the based string topology coproduct is trivial for compact symmetric spaces of higher rank and we study the implications of this result for the string topology coproduct on the free loop space.
\end{abstract}

\tableofcontents

\section{Introduction} \label{sec_introduction}

String topology is the study of algebraic structures on the homology or cohomology of the free loop space $\Lambda M$ of a closed manifold $M$.
The most prominent operation is the \textit{Chas-Sullivan product}, which is a product of the form
$$     \wedge\colon \mathrm{H}_i(\Lambda M)\otimes \mathrm{H}_j(\Lambda M)\to \mathrm{H}_{i+j-n}(\Lambda M)       $$
where $n$ is the dimension of the underlying oriented manifold $M$.
This product was first introduced by Chas and Sullivan in \cite{chas:1999} via geometric intersections of chains in the free loop space.
The geometric idea of the Chas-Sullivan product is to concatenate loops that share the same basepoint.
It turns out that the Chas-Sullivan product only depends on the oriented homotopy type of $M$, see \cite{cohen:2008}, \cite{crabb2008loop} and \cite{gruher2008generalized}.

Since the influential paper by Chas and Sullivan \cite{chas:1999} many more string topology operations have been studied.
In particular we want to mention the operation induced by the $\mathbb{S}^1$-action on the free loop space, which together with the Chas-Sullivan product induces the structure of a Batalin-Vilkovisky algebra on the homology of the free loop space, see \cite[Theorem 5.4]{chas:1999}.

In \cite{goresky:2009} Goresky and Hingston introduce the \textit{string topology coproduct}, which is a map
$$    \vee\colon \mathrm{H}_i(\Lambda M,M)\to \mathrm{H}_{i+1-n}(\Lambda M\times \Lambda M,\Lambda M\times M\cup M\times \Lambda M)\, .   $$
Here we consider $M$ as a subset of the free loop space via the identification of $M$ with the trivial loops in $M$.
The idea for this coproduct goes back to Sullivan, see \cite{sullivan:2004}.
Geometrically, the coproduct looks for loops with self-intersections at the basepoint and cuts them apart at the basepoint.
In certain situations the string topology coproduct induces a dual product in cohomology, the \textit{Goresky-Hingston product} which takes the form
$$  \ostar\colon \mathrm{H}^i(\Lambda M,M)\otimes \mathrm{H}^j(\Lambda M,M)\to \mathrm{H}^{i+j+n-1}(\Lambda M,M)\,  .  $$
Hingston and Wahl further study the string topology coproduct and the Goreky-Hingston product in \cite{hingston:2017}.
In \cite{naef2024string} Naef shows that the string topology coproduct is not a homotopy invariant in general, see also \cite{naef:2022}.
The failure of the coproduct to be homotopy invariant has recently been characterized by several authors, see \cite{naef2024simple}, \cite{kenigsberg2024obstructions} and \cite{hingston:2019}.
In particular, if the underlying manifold $M$ is simply connected, then the string topology coproduct only depends on the oriented homotopy type of $M$, see.
There is also a coproduct on the based loop space of the manifold $M$ 
$$     \vee_{\Omega} \colon \mathrm{H}_i(\Omega M,\{p_0\}) \to \mathrm{H}_{i+1-n}(\Omega M\times \Omega M,\Omega M\times \{p_0\}\cup \{p_0\}\times \Omega M)    $$
which is strongly related to the string topology coproduct.

From a geometric point of view the homology of the free loop space is very important in the study of closed geodesics in a closed Riemannian manifold $(M,g)$.
It turns out that the critical points of the energy functional $E\colon \Lambda M\to \mathbb{R}$ given by
$$   E(\gamma) = \frac{1}{2} \int_0^1 g_{\gamma(s)}(\Dot{\gamma}(s), \Dot{\gamma}(s))\,\mathrm{d}s,\quad \gamma\in\Lambda M\,,     $$
are precisely the closed geodesics in $M$ with respect to the metric $g$.
Morse-theoretic methods therefore give a connection between the topology of $\Lambda M$ and the closed geodesics in $M$.

The connection between string topology and closed geodesics has been studied in a few instances.
Hingston and Rademacher \cite{hingston:2013} use explicit computations of the Chas-Sullivan and the Goresky-Hingston product on spheres to prove a resonance property of closed geodesics.
Goresky and Hingston \cite{goresky:2009} define Chas-Sullivan and Goresky-Hingston type products in the level homologies and cohomologies of the free loop space and compute these products for spheres and projective spaces.
In particular, they then use the level products to compute the Chas-Sullivan and the Goresky-Hingston product for spheres, see \cite[Section 15]{goresky:2009}.
Moreover, Goresky and Hingston show that for spheres and projective spaces there are non-nilpotent classes for both products and that the powers of these classes corresponds to the iteration of the closed geodesics.
We also refer to \cite{hingston:2017} as well as \cite{cieliebak2020poincar} and \cite{cieliebak2023loop} for a computation of the string topology coproduct on odd-dimensional spheres.
Goresky and Hingston use the well understood geometry of a metric 
where all geodesics are closed and of the same length for their computation of the string topology products.
Note that the standard examples for such metrics are the compact symmetric spaces of rank one.

It is natural to ask how the string topology operations behave for arbitrary compact globally symmetric spaces of higher rank.
Among the higher rank symmetric spaces are the higher rank compact Lie groups, the real, complex and quaternionic Grassmannians which are not projective spaces as well as Lagrange Grassmannians and many other examples.
We refer to \cite{helgason:78} for a thorough introduction to symmetric spaces.
In this article we use explicit cycles in the based and the free loop space, respectively, of a compact symmetric space, which were introduced by Bott and Samelson \cite{bott:1958a} and by Ziller \cite{ziller:1977}.
While the Chas-Sullivan and the Goresky-Hingston product have been studied for compact symmetric spaces of rank $1$, these operations have not yet been studied extensively for compact symmetric spaces of higher rank which are not Lie groups.
In the rank $1$ case, we mentioned above the computations for spheres by Goresky and Hingston \cite{goresky:2009}, Hingston and Wahl \cite{hingston:2017} as well as by Cieliebak, Hingston and Oancea in \cite{cieliebak2020poincar} and \cite{cieliebak2023loop}.
The Chas-Sullivan product and in particular the Batalin-Vilkovisky algebra structure has been computed for spheres by Menichi in \cite{menichi2009string}.
Furthermore, the Chas-Sullivan product has been computed for complex projective spaces in \cite{cohen:2003} and for quaternionic projective spaces and the Cayley plane in \cite{cadek:2010}. 
A computation of Batalin-Vilkovisky algebra structure of complex projective spaces has been done in \cite{chataur2011loop}.
A product induced by the Chas-Sullivan product on the quotient of $\Lambda \mathbb{S}^n$ by the action generated by orientation reversal of loops has been computed by the first author in \cite{kupper:2021}.
The Goresky-Hingston product for complex and quaternionic projective spaces has been studied in \cite{stegemeyer:2022}.
The string topology of compact Lie groups has been studied by Hepworth \cite{hepworth:2009} and the second author \cite{stegemeyer:2021}.
The second author shows in \cite{stegemeyer:2021} that the string topology coproduct is trivial for compact simply connected Lie groups of rank $r\geq 2$.
Note that the Chas-Sullivan product and the string topology coproduct together are part of an infinitesimal BV bialgebra structure, see e.g. \cite{cieliebak2020poincar} and \cite{latschev2024bv}.

In the present article we study the Chas-Sullivan product for compact symmetric spaces using Ziller's explicit cycles.
Let $\gamma$ be a closed geodesic and let $\Sigma_{\gamma}\subseteq \Lambda M$ be the critical manifold containing $\gamma$.
There is a closed manifold $\Gamma_{\gamma}$, which, as we shall see later, we understand as a \textit{completing manifold} for $\Sigma_{\gamma}$.
The manifold $\Gamma_{\gamma}$ fibers over $\Sigma_{\gamma}$ and embeds into the free loop space via a map $f_{\gamma}\colon \Gamma_{\gamma}\hookrightarrow \Lambda M$.
If we take $\ZZ_2$-coefficients we take the fundamental class of $\Gamma_{\gamma}$ and we consider the class $(f_{\gamma})_*[\Gamma_{\gamma}]\in\mathrm{H}_*(\Lambda M)$.
\begin{thm}[Theorem \ref{nonnilpotent_CS_classes}]
Let $M$ be a compact symmetric space and take homology with $\ZZ_2$-coefficients. For every critical manifold $\Sigma_{\gamma}$ the associated class $\Theta =(f_{\gamma})_*[\Gamma_{\gamma}]$ is non-nilpotent in the Chas-Sullivan algebra and the powers of this class correspond to the iteration of the closed geodesics.
    More precisely, if $m_1,m_2\in \mathbb{N}$ we have
    $$   (f_{\gamma^{m_1}})_*[\Gamma_{\gamma^{m_1}}] \wedge (f_{\gamma^{m_2}})_*[\Gamma_{\gamma^{m_2}}] = (f_{\gamma^{m_1+m_2}})_*[\Gamma_{\gamma^{m_1+m_2}}]\,  .      $$
Moreover, for each $m\in\mathbb{N}$ the multiplication with $\Theta$ gives an isomorphism from the loop space homology generated by the critical manifold $\Sigma_{\gamma^m}$ to the homology generated by $\Sigma_{\gamma^{ m+1}}$.
\end{thm}
It will be made clear in Section \ref{sec_chas_sullivan} what we mean by the homology generated by the critical manifold $\Sigma_{\gamma^m}$ for $m\in\mathbb{N}$.
This theorem shows that the Chas-Sullivan product of compact symmetric spaces of higher rank is highly non-trivial.
In particular, we argue that if the rank $r$ of the compact symmetric space $M$ satisfies $r\geq 2$, then there are infinitely many critical manifolds consisting of prime closed geodesics, each of which induces a non-nilpotent homology class.
We also show that the Chas-Sullivan product of classes associated to a critical manifold $\Sigma_{\gamma^{m_1}}$ and classes associated to $\Sigma_{\gamma^{m_2}}$ is strongly related to the intersection product in the manifold $\Sigma \cong \Sigma_{\gamma^{m_1}} \cong \Sigma_{\gamma^{m_2}}$. 
See Theorem \ref{theorem_cs_intersection} for details.

Furthermore, we study the string topology coproduct on the based loop space of a compact symmetric space as well as the coproduct on the corresponding free loop space.
\begin{thm}[Corollary \ref{cor_triviality_based} and Proposition \ref{prop_wc_classes}]
    Let $M$ be a compact simply connected symmetric space of rank greater than or equal to $2$ with baseopint $p_0\in M$.
    \begin{enumerate}
        \item     The based coproduct $\vee_{\Omega}$ on the homology of the based loop space $\mathrm{H}_*(\Omega M,p_0;\ZZ_2)$ is trivial.
        \item  To every non-trivial critical manifold $\Sigma\subseteq \Lambda M$ of the energy functional one can associate a homology class $[W_{\Sigma}]\in\mathrm{H}_{\mathrm{ind}(\Sigma)}(\Lambda M,M; \ZZ_2)$, which has trivial string topology coproduct.
    \end{enumerate}
\end{thm}
Note that this result shows that the Goresky-Hingston product behaves quite differently when one compares compact symmetric spaces of rank $1$ and of higher rank, since in the rank $1$ case the dual of the class $[W_{\Sigma}]$ in cohomology is non-nilpotent in the Goresky-Hingston algebra.
In contrast to this, our results for the Chas-Sullivan product show that the Chas-Sullivan product for higher rank symmetric spaces behaves very similarly to the rank $1$ case.
It would be interesting to understand how this discrepancy between the rank $1$ case and the higher rank case fits together with the infinitesimal BV bialgebra structure discussed in \cite{cieliebak2020poincar} and \cite{latschev2024bv}.

Finally, we show that for a product of compact symmetric spaces the string topology coproduct vanishes on a large subspace of the homology of the free loop space. 

\begin{thm}[Theorem \ref{theorem_intersection_product_completing}]
    Let $M=M_1\times M_2$ be a product of two compact symmetric spaces. Let $U$ denote the subspace of $\mathrm{H}_*(\Lambda M,M;\ZZ_2)$ generated by the classes induced by completing manifolds associated to critical manifolds $\Sigma=\Sigma_1\times\Sigma_2$ where both $\Sigma_1$ and $\Sigma_2$ consist of non-constant closed geodesics in $M_1$, respectively in $M_2$. Then the string topology coproduct vanishes on $U$.
\end{thm}


This article is organized as follows.
In Section \ref{sec_def} we introduce the Chas-Sullivan product and the string topology coproduct.
The notion of a completing manifold, which is the central idea behind Bott and Samelson's and Ziller's cycles, is studied in Section \ref{sec_completing}.
In Section \ref{sec_symm_geodesics} we give an overview of geodesics in compact symmetric spaces and set the ground for the definition of the explicit cycles in the loop spaces of compact symmetric spaces in the following sections.
In Section \ref{sec_chas_sullivan} we study Ziller's cycles in detail and show how they can be used to compute the Chas-Sullivan product.
We use Bott's and Samelson's cycles in Section \ref{sec_coproduct} to study the string topology coproduct on the based and on the free loop space.
Finally, in Section \ref{sec_spheres} we study the string topology coproduct on products of compact symmetric spaces by studying how the Ziller cycles behave on product spaces.

\medskip
\noindent \textbf{Acknowledgements.} 
The authors thank the anonymous referee for their careful and thoughtful
reading of our manuscript. Their suggestions highly improved the exposition and the clarity of the article.

\section{String topology product and coproduct} \label{sec_def}

In this section we introduce the string topology operations, which we examine in this article.
We follow \cite{goresky:2009} for the definition of the Chas-Sullivan product and \cite{hingston:2017} for the definition of the string topology coproduct.

Let $(M,g)$ be a closed Riemannian manifold of dimension $n$.
We denote the unit interval by $I = [0,1]$.
We define the \textit{path space of $M$} to be $$  PM = \big\{ \gamma : I\to M\,|\, \gamma \,\,\text{absolutely continuous}, \,\, \int_0^1 g_{\gamma(t)}(\Dot{\gamma}(t),\Dot{\gamma}(t)) \,\mathrm{d}t < \infty \big\} \,   .$$
We refer to \cite[Definition 2.3.1]{klingenberg:1995} for the definition of absolute continuity of curves in a smooth manifold.
It turns out that the path space $PM$ can be given the structure of a Hilbert manifold, see \cite[Section 2.3]{klingenberg:1995}.
The \textit{free loop space of} $M$ is defined to be the subspace
$$  \Lambda M = \{\gamma \in PM \,|\, \gamma(0 ) = \gamma(1)\} ,  $$
which is in fact a submanifold of $PM$.
The underlying manifold $M$ can be seen as a submanifold of $\Lambda M$ via the identification with the trivial loops, see \cite[Proposition 1.4.6]{klingenberg:78}.
We consider the \textit{energy functional}
$$  E': PM \to [0,\infty),\qquad  E'(\gamma) =    \frac{1}{2}\int_0^1 g_{\gamma(t)}(\Dot{\gamma}(t),\Dot{\gamma}(t)) \mathrm{d}t \,     . $$
This is a differentiable function on the path space $PM$, see \cite[Theorem 2.3.20]{klingenberg:1995}.
By restriction we obtain the energy functional on the free loop space, which we will denote by $E = E'|_{\Lambda M}$.
The critical points of $E$ on $\Lambda M\setminus M$ are precisely the closed geodesics in $M$.
Moreover, for a point $p_0\in M$ we define
$$ \Omega_{p_0}M = \{ \gamma\in \Lambda M\,|\, \gamma(0) = p_0\}     $$
and call this space the \textit{based loop space of} $M$ \textit{at the point} $p_0$.
The based loop space is a submanifold of the free loop space $\Lambda M$.
Note that we will also frequently write $\Lambda$ for the free loop space $\Lambda M$ and $\Omega M$ or $\Omega$ for the based loop space if the manifold or the basepoint, respectively, are fixed.
In the following let $R$ be a commutative ring with unit and consider homology and cohomology with coefficients in $R$.
We assume that the manifold $M$ is oriented with respect to the coefficient ring $R$.
If $\epsilon>0$ is smaller than the injectivity radius of $M$, then the diagonal $\Delta M\subseteq M\times M$ has a tubular neighborhood
$$  U_M \coloneqq \{  (p,q)\in M\times M\,|\,   \mathrm{d}(p,q)<\epsilon\} \,  .     $$
Here, $\mathrm{d}$ is the distance function on $M$ induced by the Riemannian metric.
Note that the normal bundle of $M\cong \Delta M$ in the product $M\times M$ is isomorphic to the tangent bundle $TM$.
Consequently, since $U_M$ is a tubular neighborhood there is a homeomorphism of pairs 
$$   (U_M, U_M\setminus M) \cong (TM,TM\setminus M)\,  .       $$
If we pull back the Thom class of $(TM,TM\setminus M)$ via this homeomorphism we obtain a class
$$  \tau_M\in \mathrm{H}^n(U_M,U_M\setminus M) \, .      $$
On the free loop space $\Lambda M$ we consider the evaluation map $\mathrm{ev}_0\colon \Lambda\to M$ given by
$$  \mathrm{ev}_0(\gamma) = \gamma(0) \, .    $$
We define the \textit{figure-eight space} $\Lambda\times_M \Lambda$ as the pull-back
$$  \Lambda\times_M\Lambda = (\mathrm{ev}_0\times \mathrm{ev}_0)^{-1}(\Delta M) = \{(\gamma,\sigma)\in \Lambda\times \Lambda\,|\,\gamma(0) = \sigma(0)\}\,  .     $$
It turns out that the pull-back of the tubular neighborhood of the diagonal $U_M$ yields a tubular neighborhood of the figure-eight space
$$   U_{CS} \coloneqq (\mathrm{ev}_0\times\mathrm{ev}_0)^{-1}(U_M)  = \{(\gamma,\sigma)\in\Lambda\times\Lambda\,|\,\mathrm{d}(\gamma(0),\sigma(0))<\epsilon\}  \,   , $$
see e.g. \cite{hingston:2017}.
One can now pull back the Thom class $\tau_M$ via the map $\mathrm{ev}_0\times \mathrm{ev}_0$ to obtain a class
$$    \tau_{CS} \coloneqq (\mathrm{ev}_0\times\mathrm{ev}_0)^* \tau_M \in \mathrm{H}^n(U_{CS},U_{CS}\setminus \Lambda\times_M \Lambda) \, .      $$
Let $R\colon U_{CS}\to \Lambda\times_M\Lambda$ be the retraction of the tubular neighborhood.
Furthermore, for a time $s\in(0,1)$ there is a concatenation map
$$   \mathrm{concat}_s\colon \Lambda\times_M\Lambda\to \Lambda     $$
given by
$$   \mathrm{concat}_s(\gamma,\sigma)(t) \coloneqq  \begin{cases} \gamma(\tfrac{t}{s}), & 0\leq t\leq s 
\\
\sigma(\tfrac{t-s}{1-s}), & s \leq t \leq 1  \, .
\end{cases}            $$
Note that for $s,t\in(0,1)$ the maps $\mathrm{concat}_{s}$ and $\mathrm{conact}_{t}$ are homotopic via re-parametrization.
Thus, for the induced map in homology, we just write $\mathrm{concat}_*$.
We can now define the Chas-Sullivan product.
\begin{definition}
Let $M$ be a closed $R$-oriented $n$-dimensional manifold.
The \textit{Chas-Sullivan product} is defined as the composition
\begin{eqnarray*}
 \wedge\colon  \mathrm{H}_{i}(\Lambda)\otimes \mathrm{H}_j(\Lambda) &\xrightarrow{  (-1)^{n(n-i)}\times   }& \mathrm{H}_{i+j}(\Lambda\times \Lambda)
     \\ &\xrightarrow{ \hphantom{coninci}\vphantom{r}\hphantom{conci} }& \mathrm{H}_{i+j}(\Lambda\times \Lambda, \Lambda\times \Lambda\setminus \Lambda\times_M \Lambda) \\
    &\xrightarrow{\hphantom{cii}\text{excision}\hphantom{ici}}&
    \mathrm{H}_{i+j}(U_{CS}, U_{CS}\setminus \Lambda\times_M\Lambda) \\
    &\xrightarrow{\hphantom{coci} \tau_{CS}\cap \hphantom{coci} }& \mathrm{H}_{i+j-n}(U_{CS})
    \\ &\xrightarrow[]{\hphantom{conci}R_* \hphantom{conci}} & \mathrm{H}_{i+j-n}(\Lambda\times_M \Lambda) 
    \\ & \xrightarrow[]{\hphantom{ici}\mathrm{concat}_*\hphantom{ici}} & \mathrm{H}_{i+j-n}(\Lambda) \, .
    \end{eqnarray*}
\end{definition}
\begin{remark}
The Chas-Sullivan product is an associative, graded commutative and unital product, see e.g. \cite{goresky:2009}.
Moreover, it only depends on the homotopy type of the underlying manifold $M$, see \cite{cohen:2008}, \cite{crabb2008loop} and \cite{gruher2008generalized}.
\end{remark}
We now turn to the definition of the string topology coproduct.
Recall that we had the explicit tubular neighborhood of the diagonal $U_M$.
Let $\epsilon_0>0$ be a small positive number with $\epsilon_0<\epsilon$ and define 
$$  U_{M,\geq \epsilon_0} = \{(p,q)\in U_M\,|\, \mathrm{d}(p,q)\geq \epsilon_0\}\,  .    $$
Note that by \cite{hingston:2017} the pair $(U_M,U_{M,\geq\epsilon_0})$ is homeomorphic to $(TM^{<\epsilon}, TM^{<\epsilon}_{\geq\epsilon_0})$,
where 
$$   TM^{<\epsilon} \coloneqq \{ v\in TM\,|\, |v|<\epsilon\} \quad \text{and}\quad TM^{<\epsilon}_{\geq\epsilon_0} \coloneqq \{ v\in TM^{<\epsilon}\,|\, |v|\geq \epsilon_0\}      $$
and where $|\cdot|$ is the fiberwise norm induced by the Riemannian metric on $M$.
Moreover, by \cite{hingston:2017} the Thom class in $(TM,TM\setminus M)$ induces a Thom class in $(TM^{\epsilon},TM^{<\epsilon}_{\geq\epsilon_0})$.
Consequently, we obtain a class $\tau_M'\in \mathrm{H}^n(U_M,U_{M,\geq\epsilon_0})$.
Fix a basepoint $p_0\in M$ and consider the sets
$$   B_{p_0} \coloneqq \{q\in M\,|\, \mathrm{d}(p_0,q)<\epsilon\}\quad \text{and}\quad B_{p_0,\geq\epsilon_0} \coloneqq \{ q\in B_{p_0}\,|\, \mathrm{d}(p_0,q)\geq \epsilon_0\} \, .$$
Let $i$ be the inclusion
$$      i\colon B_{p_0} \xrightarrow[]{\cong} \{p_0\}\times B_{p_0} \hookrightarrow  U_M\,  . $$
This induces a map of pairs $j\colon (B_{p_0},B_{p_0,\geq\epsilon_0})\to (U_M,U_{M,\geq\epsilon_0})$.
If we pull back the class $\tau_M'$ via this map we obtain a class
$$   \tau_{p_0} \coloneqq j^*\tau_M' \in \mathrm{H}^n(B_{p_0},B_{p_0,\geq\epsilon_0}) \cong \mathrm{H}(M,M\setminus\{p_0\})      $$
and by the properties of a Thom class this class coincides with the generator of $\mathrm{H}(M,M\setminus\{p_0\})$ induced by the orientation of $M$.

On the product of the free loop space with the unit interval $\Lambda\times I$ there is an evaluation map
$$  \mathrm{ev}_I\colon \Lambda\times I\to M\times M     $$
given by $$   \mathrm{ev}_I(\gamma,s)  = (\gamma(0),\gamma(s)) \, .   $$
Define
\begin{equation*}
  U_{GH} = (\mathrm{ev}_I)^{-1}(U_M) = \{ (\gamma,s) \in \Lambda M\times I\,|\, \mathrm{d}(\gamma(0),\gamma(s))< \epsilon\} 
  \end{equation*}
  and
  \begin{equation*}
  U_{GH,\geq\epsilon_0} =  (\mathrm{ev}_I)^{-1}(U_{M,\geq\epsilon_0}) =
  \{ (\gamma,s)\in U_{GH}\,|\, \mathrm{d}(\gamma(0),\gamma(s))\geq \epsilon_0\}\,  .
\end{equation*}
The set $U_{GH}$ is an open neighborhood of the subspace 
$$ F_{GH} \coloneqq  (\mathrm{ev}_I)^{-1}(\Delta M)  =   \{ (\gamma,s)\in \Lambda\times I\,|\, \gamma(0)=\gamma(s)\}\,   . $$
The evaluation $\mathrm{ev}_I$ defines a map of pairs 
$$ \mathrm{ev}_I\colon    (U_{GH} ,U_{GH,\geq\epsilon_0}) \to (U_M,U_{M,\geq\epsilon_0})     $$
and we can pull back the class $\tau_M'$ to obtain a class
$$     \tau_{GH} \coloneqq (\mathrm{ev}_I)^* \tau_M'\in \mathrm{H}^n( U_{GH} ,U_{GH,\geq\epsilon_0})\,  .       $$

The above constructions can be done similarly for the based loop space, i.e. we obtain an open neighborhood $ U^{\Omega}_{GH} = U_{GH}\cap (\Omega\times I) $ of the space
$$  F^{\Omega}_{GH}  \coloneqq F_{GH}\cap (\Omega\times I) = \{ (\gamma,s)\in \Omega\times I\,|\, \gamma(s) = p_0\}      $$
and we pull back the class $\tau_{p_0}$ to a class
$$  \tau^{\Omega}_{GH} \in \mathrm{H}^n(U^{\Omega}_{GH},U^{\Omega}_{GH,\geq\epsilon_0})     $$
where $U^{\Omega}_{GH,\geq\epsilon_0} \coloneqq U_{GH,\geq\epsilon_0} \cap (\Omega\times I)$.

As Hingston and Wahl argue in \cite{hingston:2017}, there is a retraction map $R_{GH}\colon U_{GH}\to F_{GH}$ and it is easy to see that this restricts to a retraction $U_{GH}^{\Omega}\to F^{\Omega}_{GH}$.
Furthermore, we have a cutting map $\mathrm{cut}\colon F_{GH}\to \Lambda\times \Lambda$ given by
$$   \mathrm{cut}(\gamma,s) = (\gamma|_{[0,s]},\gamma|_{[s,1]}) ,    $$
where by $\gamma|_{[0,s]}$ we mean a reparametrized version of the restriction $\gamma|_{[0,s]}$.
Again, this restricts appropriately to a based version.

Let $[I]$ be the positively oriented generator of $\mathrm{H}_1(I,\partial I)$ with respect to the standard orientation of the unit interval.
We write $p_0$ for the basepoint of the based loop space which is just the trivial loop at the point $p_0$ and by abuse of notation we shall denote the set consisting of the single element $p_0$ just by $p_0$. 
\begin{definition} Let $M$ be a closed $R$-oriented manifold of dimension $n$.
\begin{enumerate}
    \item 
The \textit{string topology coproduct} is defined as the map
\begin{eqnarray*}
 \vee: \mathrm{H}_{*}(\Lambda,M) &\xrightarrow{\hphantom{oi}\times [I]\hphantom{oi}}& \mathrm{H}_{* +1}(\Lambda\times I, \Lambda\times\partial I\cup M\times I)
     \\ &\xrightarrow{ \hphantom{o} \tau_{GH}\cap \hphantom{o} }& \mathrm{H}_{* +1-n}(U_{GH}, \Lambda\times\partial I\cup p_0\times I) 
     \\
    &\xrightarrow{ (\mathrm{R}_{GH})_* }&
    \mathrm{H}_{* +1-n}(F_{GH}, \Lambda\times\partial I\cup M\times I) \\
    &\xrightarrow{ \hphantom{i} (\mathrm{cut})_* \hphantom{i} }& \mathrm{H}_{* +1-n}(\Lambda\times\Lambda, \Lambda\times M\cup M\times \Lambda)\, . 
    \end{eqnarray*}
    \item
    The \textit{based string topology coproduct} is defined as the map
    \begin{eqnarray*}
 \vee_{\Omega}: \mathrm{H}_{*}(\Omega,p_0) &\xrightarrow{ \hphantom{oi} \times [I] \hphantom{oi} }& \mathrm{H}_{* +1}(\Omega\times I, \Omega\times\partial I\cup p_0\times I)
     \\
     &\xrightarrow{\hphantom{o} \tau_{GH}^{\Omega}\cap \hphantom{o} }& \mathrm{H}_{* +1-n}(U_{GH}^{\Omega}, \Omega\times\partial I\cup p_0\times I) \\
    &\xrightarrow{(\mathrm{R}_{GH})_*}&
    \mathrm{H}_{* +1-n}(F^{\Omega}_{GH}, \Omega\times\partial I\cup p_0\times I) \\
    &\xrightarrow{\hphantom{i} (\mathrm{cut})_* \hphantom{i}}& \mathrm{H}_{* +1-n}(\Omega\times\Omega, \Omega\times p_0\cup p_0\times \Omega)\, . 
    \end{eqnarray*}
    
\end{enumerate}
\end{definition}
\begin{remark}
If we take $R=\mathbb{F}$ to be a field, then by the Künneth isomorphism the string topology coproduct induces a map
$$   \vee' \colon \mathrm{H}_i(\Lambda,M) \to (\mathrm{H}_{*}(\Lambda,M)\otimes \mathrm{H}_{*}(\Lambda,M))_{i+1-n} \, .       $$
After a sign correction this map is indeed coassociative and graded cocommutative, see \cite[Theorem 2.14]{hingston:2017} and therefore justifies the name \enquote{coproduct}.
\end{remark}

Note that the based string topology coproduct and the string topology coproduct on the free loop space are compatible with respect to the inclusion $i\colon \Omega \hookrightarrow \Lambda$.
More precisely, we have the commutativity of the following diagram, see \cite[Proposition 2.4]{stegemeyer:2021}.
$$
\begin{tikzcd}    \mathrm{H}_i(\Omega,p_0) \arrow[]{r}{i_*} \arrow[]{d}{\vee_{\Omega}} & \mathrm{H}_i(\Lambda ,M) \arrow[]{d}{\vee}
    \\
    \mathrm{H}_{i+1-n}(\Omega \times \Omega,\Omega\times p_0\cup p_0\times \Omega) \arrow[]{r}{(i\times i)_*} & \mathrm{H}_{i+1-n}(\Lambda \times \Lambda,\Lambda\times M\cup M\times \Lambda) .
\end{tikzcd}
$$
By dualizing, the coproduct with coefficients in a field induces a product in cohomology.
\begin{definition}
Let $\mathbb{F}$ be a field and take homology and cohomology with coefficients in $\mathbb{F}$.
Let $M$ be a closed $\mathbb{F}$-oriented manifold.
If $\alpha,\beta\in \mathrm{H}^{*}(\Lambda M,M;\mathbb{F})$ are relative cohomology classes, then the \textit{Goresky-Hingston product} $\alpha\ostar \beta$ is defined by taking the dual to the string topology coproduct, i.e. 
$$    \alpha \ostar \beta : = \vee^* (\alpha\otimes \beta) = (\alpha\otimes \beta) \circ \vee.    $$
Here, we make the canonical identification $\mathrm{H}^{*}(\Lambda M,M;\mathbb{F}) \cong \mathrm{Hom}(\mathrm{H}_{*}(\Lambda M,M;\mathbb{F}),\mathbb{F})$.
\end{definition}

Note that the Goresky-Hingston product can also be defined for more general coefficients, see \cite[Section 1.6]{hingston:2017}.

We conclude this section by defining the basepoint intersection multiplicity of a relative homology class for both the based loop space and the free loop space.
Hingston and Wahl define the notion of \textit{basepoint intersection multiplicity} in \cite[Section 5]{hingston:2017} only in the case of the free loop space.
The analogous concept for the based loop space is straight-forward.

\begin{definition} Let $M$ be a closed manifold with basepoint $p_0\in M$.
\begin{enumerate}
    \item 
        Let $[X]\in \mathrm{H}_{*}(\Lambda M,M; R)$ be a homology class with representing cycle $X\in \mathrm{C}_{*}(\Lambda M,M)$.
         Assume that the relative cycle $X$ is itself represented by a cycle $x\in \mathrm{C}_{*}(\Lambda M)$.
        The \textit{basepoint intersection multiplicity} $\mathrm{int}([X])$ of the class $[X]$ is the number
     $$    \mathrm{int}([X]) \coloneqq \inf_{A \sim x} \Big(\sup\big[ \#(\gamma^{-1}(\{ \gamma(0) \})) \,|\, \gamma\in \mathrm{im}(A),\,\,\, E(\gamma)> 0 \big] \Big)  -1   $$
    where the infimum is taken over all cycles $A \in \mathrm{C}_{*}(\Lambda M)$ homologous to $x$.
    \item
     Let $[X]\in \mathrm{H}_{*}(\Omega_{p_0} M,p_0; R)$ be a homology class with representing cycle $X\in \mathrm{C}_{*}(\Omega_{p_0} M,p_0)$.
     Assume that the relative cycle $X$ is itself represented by a cycle $x\in \mathrm{C}_{*}(\Omega_{p_0} M)$.
    The \textit{basepoint intersection multiplicity} $\mathrm{int}([X])$ of the class $[X]$ is the number
    $$    \mathrm{int}([X]) \coloneqq \inf_{A \sim x} \Big(\sup\big[ \#(\gamma^{-1}(\{p_0\})) \,|\, \gamma\in \mathrm{im}(A),\,\,\, E(\gamma)> 0 \big] \Big)  -1   $$
    where the infimum is taken over all cycles $A \in \mathrm{C}_{*}(\Omega_{p_0} M)$ homologous to $x$.
\end{enumerate}
\end{definition}

Note that in the present paper loops are maps $\gamma\colon [0,1]\to M$ and therefore any loop has a basepoint intersection at $t=0$ and $t= 1$.
We therefore subtract $1$ in the above definition in order not to count this trivial intersection twice.
In \cite{hingston:2017} loops are considered as maps $\mathbb{S}^1\to M$ and therefore the additional term $-1$ does not show up in \cite{hingston:2017}.

By considering finite-dimensional models of $\Lambda M$, respectively of $\Omega_{p_0}M$, one can see that the intersection multiplicity of a homology class is finite.
Hingston and Wahl prove in the case of the free loop space that if a homology class $Y\in \mathrm{H}_{*}(\Lambda M,M)$ has intersection multiplicity $1$ then its coproduct vanishes, see
\cite[Theorem 3.10]{hingston:2017}.
Since the proof of this theorem which includes Lemma 3.5 and Proposition 3.6 of \cite{hingston:2017} only uses maps that keep the basepoint of the loops fixed, the proof can be transferred directly to the case of the based loop space.
We obtain the following result.
\begin{prop}[\cite{hingston:2017}, Theorem 3.10] \label{prop_intersection_trivial}
Let $M$ be a closed manifold with basepoint $p_0\in M$ and let $R$ be a commutative ring.
Assume that $M$ is $R$-oriented.
\begin{enumerate}
    \item Suppose that $[X]\in \mathrm{H}_{*}(\Lambda M,M;R)$ has basepoint intersection multiplicity $\mathrm{int}([X]) = 1$, then the string topology coproduct $\cpr[X]$ vanishes.
    \item Suppose that $[X]\in \mathrm{H}_{*}(\Omega_{p_0}M,p_0;R)$ has basepoint intersection multiplicity $\mathrm{int}([X]) = 1$, then the based coproduct $\cpr_{\Omega}[X]$ vanishes.
\end{enumerate}
\end{prop}

\section{Completing manifolds} \label{sec_completing}

In this section we introduce the concept of a completing manifold.
This will be the key technique in our study of the string topology operations on symmetric spaces in Sections \ref{sec_coproduct} and \ref{sec_chas_sullivan}.
We closely follow \cite{hingston:2013oancea} and \cite{oancea:2015}.

Assume that $X$ is a Hilbert manifold and that $F \colon X\to \RR$ is a Morse-Bott function.
Let $a$ be a critical value of $F$ and assume that the set of critical points at level $a$ is a non-degenerate submanifold $\Sigma$ with index $\lambda$ and with orientable negative bundle.
As usual in Morse theory we denote by $X^{\leq a}$ the set
$$    X^{\leq a} \coloneqq \{ x\in X\,|\, F(x)\leq a\}\,  .       $$
Similarly, one defines the strict sublevel set $X^{<a}$.
It is well-known that the level homology $\mathrm{H}_{*}(X^{\leq a},X^{<a})$ is completely determined by the homology of the critical manifold $\Sigma$.
More precisely, we have
$$     \mathrm{H}_{*}(X^{\leq a},X^{<a}) \cong \mathrm{H}_{* - \lambda}(\Sigma) \, .      $$
In general, there is no direct link between the homology of $\Sigma$ and the homology of the sublevel set $X^{\leq a}$.
As we shall see in the following the existence of a completing manifold however is sufficient for the property that the homology of $X^{\leq a}$ is a direct sum of the homology of the strict sublevel set $X^{<a}$ and the homology of the critical manifold $\Sigma$. 
\begin{definition} \label{def_completing_mfld}
Let $X$ be a Hilbert manifold and let $F\colon X\to \RR$ be a Morse-Bott function.
Suppose that $a$ is a critical value of $F$ and assume that $\Sigma$ is the orientable non-degenerate critical submanifold at level $a$ of index $\lambda$ and of dimension $m = \mathrm{dim}(\Sigma)$.
A \textit{completing manifold} for $\Sigma$ is a closed, orientable manifold $\Gamma$ of dimension $\lambda + m$
with a closed submanifold $s\colon S\hookrightarrow \Gamma$ of codimension $\lambda$ and a map $f: \Gamma \to X^{\leq a}$ such that the following holds.
The map $f$ is an embedding near $S$ with the restriction $f|_S$ mapping $S$ homeomorphically onto $\Sigma$ and such that $f^{-1}(\Sigma) = S$.
Furthermore, there is a retraction map $p:\Gamma\to S$, i.e. $p \circ s = \mathrm{id}_S$, and the map $f$ induces a map of pairs
$$  f \colon (\Gamma,\Gamma\setminus S) \to (X^{\leq a}, X^{<a}) \, .    $$
In case that we consider homology with $\ZZ_2$-coefficients we drop the orientability conditions on $\Sigma$ and $\Gamma$.
\end{definition}
\begin{remark}
Note that the above definition is the one of a \textit{strong} completing manifold in \cite{oancea:2015}.
Since in the present article we only encounter this stronger version of a completing manifold, we just call it \textit{completing manifold}.
Furthermore, we remark that in the present article all maps $f\colon \Gamma\to X$ from the completing manifold $\Gamma$ to the Hilbert manifold $X$ will be embeddings so that the property of being an embedding near $S$ is clearly satisfied.
\end{remark}
Let $g \colon M\to N$ be a map between closed oriented manifolds.
We recall that the Gysin map
$$     g_!\colon \mathrm{H}_i(N) \to \mathrm{H}_{i+\mathrm{dim}(M)-\mathrm{dim}(N)}(M)$$
is defined by
$$g_! \colon \mathrm{H }_i(N) \xrightarrow[]{(PD_{N})^{-1}}  \mathrm{H}^{\mathrm{dim}(N) - i}(N) \xrightarrow[]{g^*} \mathrm{H}^{\mathrm{dim}(N)-i}(M) \xrightarrow[]{PD_{M}} \mathrm{H}_{\mathrm{dim}(M) - (\mathrm{dim}(N)-i)}(M) , $$
where 
$$\mathrm{PD}_L\colon \mathrm{H}^j(L) \xrightarrow[]{\cong} \mathrm{H}_{\mathrm{dim}(L)-j}(L)  $$ stands for Poincaré duality on the closed oriented manifold $L$.
In particular, we note that Gysin maps behave contravariantly, i.e. if $h\colon N\to L$ is a second map between compact oriented manifolds, then
\begin{equation} \label{eq_functoriality_gysin}
        (h\circ g)_! =   g_!\circ h_! \, .      
\end{equation}
If we use $\ZZ_2$-coefficients then all orientability conditions in the above can be dropped.

Let $X$ be a Hilbert manifold with a Morse-Bott function $F\colon X\to \RR$.
Assume that $\Sigma$ is an oriented critical submanifold at level $a$ and that $\Gamma$ is an oriented completing manifold for $\Sigma$.
Then we have the embedded submanifold $s\colon S\hookrightarrow \Gamma$ and by definition there is retraction map $p\colon \Gamma\to S$, i.e. we have
\begin{equation} \label{eq_retraction_inclusion}
      p\circ s = \id_S\,  .    
\end{equation}
Moreover, note that we get a map
$$     \mathrm{H}_i(\Gamma)\to \mathrm{H}_i(\Gamma,\Gamma\setminus S)\xrightarrow[]{\cong} \mathrm{H}_{i-\lambda}(S)     $$
where the first map is induced by the inclusion of pairs and the second one is the composition of excision and the Thom isomorphism.
It turns out, see \cite[Theorem 11.3]{bredon:2013}, that this map coincides up to sign with the Gysin map $s_!$.
Consequently, by equations \eqref{eq_functoriality_gysin} and \eqref{eq_retraction_inclusion} we see that $p_!$ is a right inverse for this map up to sign.
In particular, the long exact sequence of the pair $(\Gamma,\Gamma\setminus S)$ yields a short exact sequence in each degree, which splits via the map $p_!$.
Thus, for each $i\in \NN$ we have a short exact sequence
$$   0 \to \mathrm{H}_{i}(\Gamma\setminus S) \to \mathrm{H}_i(\Gamma) \xrightarrow[]{s_!} \mathrm{H}_{i-\lambda}(S) \to 0     $$
and due to the splitting $p_!$ we see that
$$   \mathrm{H}_i(\Gamma) \cong \mathrm{H}_{i}(\Gamma\setminus S)\oplus \mathrm{H}_{i-\lambda}(S) \cong \mathrm{ker}(s_!)\oplus \mathrm{H}_{i-\lambda}(S)\cong\mathrm{ker}(s_!)\oplus \mathrm{im}(p_!) \, .     $$
Combining these ideas with the existence of the map $f\colon \Gamma\to X$ one can show the following.

\begin{prop}[\cite{oancea:2015}, Lemma 6.2] \label{prop_completing_mfld}
Let $X$ be a Hilbert manifold, $F\colon X\to \RR$ a Morse-Bott function and let $a$ be a critical value of $F$.
Assume that the set of critical points at level $a$ is an oriented non-degenerate critical submanifold $\Sigma$ of index $\lambda$.
If there is a completing manifold for $\Sigma$ then
\begin{equation} \label{eq_splitting_completing_manifold}
       \mathrm{H}_{*}(X^{\leq a}) \cong \mathrm{H}_{*}(X^{<a}) \oplus \mathrm{H}_{* }(X^{\leq a},X^{<a}) \cong \mathrm{H}_{*}(X^{<a}) \oplus \mathrm{H}_{* -\lambda}(\Sigma) \,  .         
\end{equation}
\end{prop}
\begin{remark}\label{splitting_sublevels}
Note that the construction of a completing manifold does not only give us the isomorphism of Proposition \ref{prop_completing_mfld}, but it also yields an explicit way of describing the homology classes which \textit{come from} the critical manifold $\Sigma$.
More precisely, by the above discussion it is clear that the composition
$$     f_*\circ p_! \colon \mathrm{H}_i(\Sigma)\to \mathrm{H}_{i+\lambda}(\Gamma)\to \mathrm{H}_{i+\lambda}(X^{\leq a})    $$
is injective.
Here, the map $f$ is the map of the completing manifold $\Gamma$ into $X$ and $p$ is the retraction $\Gamma\to S\cong \Sigma$.

Recall that the homology of the completing manifold $\Gamma$ splits as 
$$ \mathrm{H}_i(\Gamma)\cong \mathrm{H}_i(\Gamma\setminus S)\oplus \mathrm{H}_{i}(\Gamma,\Gamma\setminus S)  \cong \mathrm{H}_i(\Gamma\setminus S)\oplus \mathrm{H}_{i-\lambda}(\Sigma)\,   .  $$
We note that, while we usually only care about the classes in $\mathrm{H}_i(\Gamma)$ that are in the image of $p_!\colon \mathrm{H}_i(\Sigma)\to \mathrm{H}_{i+\lambda}(\Gamma)$, we cannot say much about the image of $\mathrm{H}_i(\Gamma\setminus S)$ under the map $f_*$.
For example, it is not clear in general when a class in $\mathrm{H}_i(\Gamma\setminus S)$ lies in the kernel of $f_*$.
However, by the commutativity of the diagram
\begin{equation*}
    \xymatrix{
    \mathrm{H}_i(X^{<a})  \ar[r]&\mathrm{H}_i(X^{\leq a}) \ar[r] &  \mathrm{H}_i(X^{\leq a},X^{<a})   \\
    \mathrm{H}_i(\Gamma\setminus S) \ar[u]^-{f_*} \ar[r] &  \mathrm{H}_i(\Gamma) \ar[u]^-{f_*} \ar[r] & \mathrm{H}_i(\Gamma,\Gamma\setminus S) \ar[u]^-{f_*}
    }
\end{equation*}
we do see that $\mathrm{H}_i(\Gamma\setminus S)$ is mapped to the homology of the strict sublevel set $\mathrm{H}_i(X^{<a})$ via $f_*$.
\end{remark}

If $X$ is a Hilbert manifold with a Morse-Bott function $F\colon X\to \mathbb{R}$ let us assume that for every critical value $a$ we have
\begin{equation} \label{eq_perfect_morse}
        \mathrm{H}_{*}(X^{\leq a}) \cong \mathrm{H}_{*}(X^{<a}) \oplus \mathrm{H}_{*}(X^{\leq a},X^{<a})\,        .
\end{equation}
In this case we say that the function $F$ is a \textit{perfect} Morse-Bott function.
This property clearly holds if every critical submanifold is oriented and has a completing manifold.
Ziller shows in \cite{ziller:1977} that the energy functional on the free loop space of a compact symmetric space is perfect if one takes $\ZZ_2$-coefficients.
He uses explicit completing manifolds, which we shall study in detail in Section \ref{sec_chas_sullivan} in order to compute the Chas-Sullivan product on symmetric spaces partially.

\section{Geodesics in symmetric spaces} \label{sec_symm_geodesics}

In this section we review some basic constructions in the theory of globally symmetric spaces.
In Sections \ref{sec_chas_sullivan} and \ref{sec_coproduct} we shall use the concepts from this section and therefore want to introduce the necessary tools in a concise manner.
A general reference for the geometry of symmetric spaces is \cite{helgason:78}.

A connected Riemannian manifold $(M,g)$ is a \textit{symmetric space} if for each point $p\in M$ there is an isometry $s_p\colon M\to M$ that fixes $p$ and such that $(\mathrm{D}s_p)_p = -\mathrm{id}_{T_p M}$.
The isometry $s_p$ is also called \emph{geodesic involution} at $p$.
Note that the isometry group of any closed Riemannian manifold is a compact Lie group.
Let $G$ be the connected component of the identity of the group of isometries of the symmetric space $(M,g)$ and fix a basepoint $o\in M$.
There is an involutive Lie group automorphism $S\colon G\to G$ which is given by $S(\varphi) = s_o\circ \varphi\circ s_o$.
Since we assume that $M$ is connected, the group $G$ acts transitively on $M$, see \cite[Theorem IV.3.3]{helgason:78}. 
Hence, $M$ is a homogeneous space $M = G/K$ with $K$ being the stabilizer of $o$ under the action of $G$.
It further holds that 
$$  \mathrm{Fix}(S)_0 \subseteq K\subseteq \mathrm{Fix}(S) \, .   $$
Here, $\mathrm{Fix}(S)$ denotes the subgroup of $G$ consisting of elements which are fixed by the involutive automorphism $S$, i.e.
$$   \mathrm{Fix}(S) = \{k\in G\,|\,S(k) = k\} \, .    $$
Moreover, by the subscript $0$ we denote the connected component of the unit element $e\in G$.

Let $\mathfrak{g} \cong T_e G$ be the Lie algebra of $G$.
The differential $\mathrm{D}S_e\colon \mathfrak{g}\to \mathfrak{g}$ is an involutive Lie algebra automorphism of $\mathfrak{g}$.
Therefore, $\mathfrak{g}$ has a decomposition of the form
$$  \mathfrak{g} = \mathfrak{k}\oplus \mathfrak{m}\, ,   $$
with 
$$  \mathfrak{k} = \mathrm{Eig}(\mathrm{D}S_e\,\,,+1) \quad \text{and}\quad \mathfrak{m} = \mathrm{Eig}(\mathrm{D}S_e\,\,,-1)\,     $$
where $\mathrm{Eig}(DS_e,\pm 1)$ is the $\pm 1$-eigenspace in $\mathfrak{g}$ of the linear map $\mathrm{D}S_e$.
The subspace $\mathfrak{k}$ is precisely the Lie algebra of the subgroup $K$ and sits inside $\mathfrak{g}$ as a Lie subalgebra, see \cite[Theorem IV.3.3]{helgason:78}.
Moreover, the Lie bracket of $\mathfrak{g}$ restricted to $\mathfrak{m}$ satisfies $[\mathfrak{m},\mathfrak{m}]\subseteq \mathfrak{k}$.
The adjoint action $\mathrm{Ad}\colon G\times \mathfrak{g}\to \mathfrak{g}$ of $G$ restricts to an action of $K$ on $\mathfrak{m}$.
To see this, let $k\in K$ and $g\in G$.
Since $K\subseteq \mathrm{Fix}(S)$ we have that
$$   S (kgk^{-1}) = s_o kgk^{-1}  s_o = s_o k s_o s_o g s_o s_o k^{-1} s_o   =  S(k)  s_o g s_o  S(k^{-1}) = k  S(g)  k^{-1} .   $$
Hence, it holds that $S\circ \mathrm{Conj}_k = \mathrm{Conj}_k \circ S$ and by differentiating this equation at $e\in G$ one checks that the eigenspaces of $\mathrm{D}S_e$ are preserved by $(\mathrm{D}\mathrm{Conj}_k )_e = \mathrm{Ad}_k\colon \mathfrak{g}\to \mathfrak{g}$.
Note that if $\pi\colon G\to M = G/K$ is the canonical projection then
$$   \mathrm{D}\pi_e |_{\mathfrak{m}} \colon \mathfrak{m} \to T_o M   $$
is an isomorphism, see \cite[Theorem IV.3.3]{helgason:78}.
Moreover, the action of $K$ on $T_o M$ which is induced by differentiating the action of $K$ on $M$ is equivalent to the adjoint action of $K$ on $\mathfrak{m}$ via $\mathrm{D}\pi_e|_{\mathfrak{m}}$.

Next, we discuss the roots of a compact symmetric space.
Since the structure theory of Lie algebras behaves particularly well for complex Lie algebras we consider the complexification $\mathfrak{g}_{\mathbb{C}}$ of $\mathfrak{g}$ for a moment.
One introduces the notion of a \textit{Cartan subalgebra} $\mathfrak{h}\subseteq \mathfrak{g}_{\mathbb{C}}$ of $\mathfrak{g}_{\mathbb{C}}$, see \cite[p. 162]{helgason:78} for a definition.
This is in particular a maximal abelian subalgebra of $\mathfrak{g}_{\mathbb{C}}$.
Note that the adjoint maps $\mathrm{ad}_{H_1},\mathrm{ad}_{H_2}\colon \mathfrak{g}\to \mathfrak{g}$ commute for $H_1,H_2\in\mathfrak{h}$ since $\mathfrak{h}$ is abelian and thus $$\mathrm{ad}_{H_1}\circ \mathrm{ad}_{H_2} - \mathrm{ad}_{H_2}\circ \mathrm{ad}_{H_1} = \mathrm{ad}_{[H_1,H_2]} = 0 .$$
Consequently, the maps $\mathrm{ad}_H\colon \mathfrak{g}\to \mathfrak{g}$, $H\in\mathfrak{h}$ can be diagonalized simultaneously.
We want to control the eigenvalues of $\mathrm{ad}_H$ as $H$ varies in $\mathfrak{h}$ which leads to the notion of a root.
A \textit{root of }$\mathfrak{g}_{\mathbb{C}}$ is a non-zero element $\alpha$ of the dual space $\mathfrak{h}^*$ such that there is a non-trivial element $X\in\mathfrak{g}_{\mathbb{C}}$ with
$$   \mathrm{ad}_H(X) =  [H,X] =\pi \, i \, \alpha(H) X \quad \text{for all} \,\,\, H\in \mathfrak{h}\, .     $$
We denote the set of roots of $\mathfrak{g}_{\mathbb{C}}$ by $\Delta'$.
By definition the \emph{root subspace}
$$  \mathfrak{g}^{\alpha} = \{ X\in \mathfrak{g}_{\mathbb{C}} \,|\, [H,X] = \pi \,i\, \alpha(H) X \}    $$
is a non-trivial subspace of $\mathfrak{g}_{\mathbb{C}}$.
It turns out that if $\mathfrak{g}_{\mathbb{C}}$ semisimple, then $\mathfrak{g}_{\mathbb{C}}$ decomposes as the direct sum of the Cartan subalgebra $\mathfrak{h}$ and the direct sum over all root subspaces $\mathfrak{g}^{\alpha}$, $\alpha\in \Delta'$, i.e.
$$  \mathfrak{g}_{\mathbb{C}} = \mathfrak{h} \oplus \bigoplus_{\alpha\in\Delta'} \mathfrak{g}^{\alpha} ,   $$
see \cite[Theorem III.4.2]{helgason:78}.
The set of roots therefore has to be finite.

Coming back to the geometry of the symmetric space $(M,g)$ note that the subspace $\mathfrak{m}\subseteq \mathfrak{g}$ which we had identified with the tangent space $T_o M$ via $\mathrm{D}\pi_e$ can also be seen as a subspace of the complexified Lie algebra $\mathfrak{m}\subseteq \mathfrak{g}\subseteq \mathfrak{g}_{\mathbb{C}}$.
It turns out that there is a maximal abelian subspace $\mathfrak{a} \subseteq \mathfrak{m}$ of $\mathfrak{m}$ such that $\mathfrak{a}$ is contained in a Cartan subalgebra $\mathfrak{h}$ of $\mathfrak{g}_{\mathbb{C}}$, see \cite[p. 284]{helgason:78}.
Recall that $\mathfrak{m}$ is not a Lie algebra, but rather $[\mathfrak{m},\mathfrak{m}]\subseteq \mathfrak{k}$.
By an \emph{abelian subspace} of $\mathfrak{m}$ we mean in this case that $\mathfrak{a}$ satisfies $[\mathfrak{a},\mathfrak{a}]= \{0\}$.
We call the subspace $\mathfrak{a}$ a \textit{Cartan subalgebra} of $(G,K)$.
Let $\Delta$ be the set of roots $\alpha\in \Delta'$ such that $\alpha|_{\mathfrak{a}} \neq 0$ where we identify two roots $\alpha_1,\alpha_2\in \Delta'$ if $\alpha_1|_{\mathfrak{a}} = \alpha_2|_{\mathfrak{a}}$.
We say that $\Delta$ is the set of \emph{roots} of $(M,g)$. 
We remark that the roots of the symmetric space $(M,g)$ are real-valued and they govern a lot of the geometric features of $M$.
In particular the conjugate points along geodesics can be read off from the root system as we shall see later.

Further, we define the map $\mathrm{Exp}\colon \mathfrak{m}\to M$ by
$$  \mathrm{Exp} =   \pi\circ \mathrm{exp}|_{\mathfrak{m}}      $$
with $\mathrm{exp}$ being the Lie group exponential of the Lie group $G$.
It turns out that this map is in fact the Riemannian exponential map at the point $p$ under the canonical identification $T_o M\cong \mathfrak{m}$, see \cite[Theorem IV.3.3]{helgason:78}.
Set $T = \mathrm{Exp}(\mathfrak{a})\subseteq M$.
This is a flat totally geodesic and homogeneous submanifold and therefore isometric to an embedded flat torus.
We call $T$ the \textit{maximal torus}.
In particular, the geodesics in $T$ are just the images of the straight lines in $\mathfrak{a}$ under the exponential map $\mathrm{Exp}$.

Before we move on to the behavior of geodesics in symmetric spaces, we want to introduce the notion of the rank of a symmetric space.
\begin{definition}
    Let $M$ be a compact symmetric space.
    The dimension of the maximal abelian subspace $\mathfrak{a}\subseteq \mathfrak{m}$ as introduced above is called the \emph{rank} of $M$.
\end{definition}
Note that the rank of $M$ does not depend on the choice of $\mathfrak{a}$, see \cite[Lemma V.6.3]{helgason:78}.
The compact symmetric spaces of rank $1$ are precisely the spheres and the real, complex and quaternionic projective spaces as well as the Cayley plane.
Further examples of compact symmetric spaces are compact Lie groups; real, complex and quaternionic Grassmannians as well as Lagrange Grassmannians and the space of orthogonal complex structures in $\mathbb{R}^{2n}$.
In fact symmetric spaces can be classified by extending the classification results for simple Lie algebras to the setting of symmetric spaces.
See \cite[Chapter X]{helgason:78} for details.

\begin{figure}[t]
\centering
\includegraphics[scale=0.14]{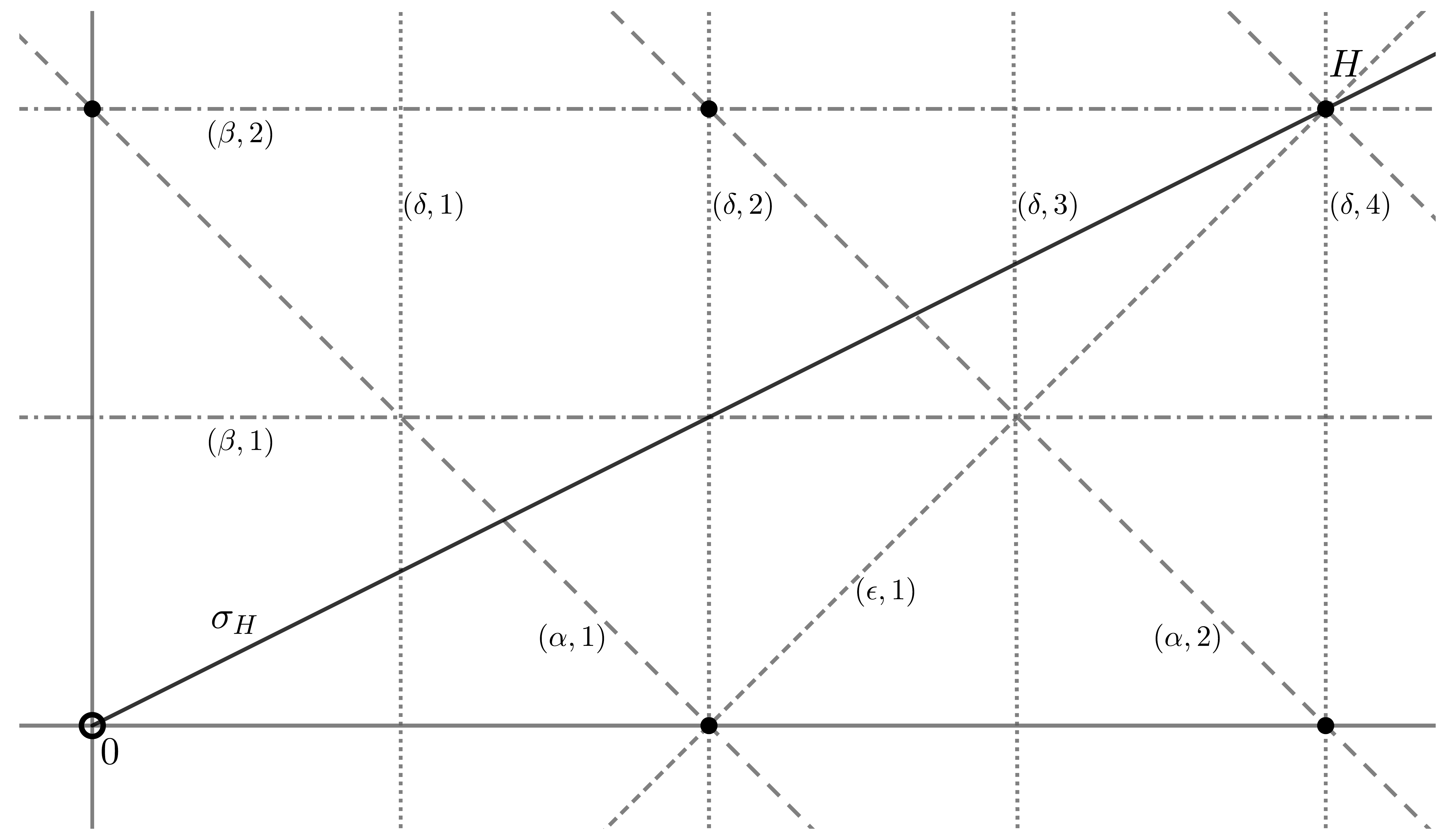}
\caption{Example of the maximal abelian subspace of a symmetric space of rank $2$.
The above figure shows the maximal abelian subspace for the complex Grassmannian $\mathrm{Gr}_2(\CC^4)$, see \cite[Section 4.2]{Sakai3} and \cite[Example 5.5]{Sakai1}.
In this example there are four positive roots 
$\{\alpha,\beta,\delta,\epsilon\}$ and the corresponding singular planes are drawn as dashed or dotted lines. 
The lattice points are pictured as dots.
The ray $\sigma_H$ is mapped to a closed geodesic by $\Exp$ since it intersects the lattice $\mathcal{F}$ at the point $H$. Observe that there are five conjugate points in the interior of the corresponding closed geodesic $\gamma_H$. They can be read off by considering the intersections of $\sigma_H$ with the singular planes.}
\label{FigureRootSystem}
\end{figure}

We now want to see what the roots of the compact symmetric space $(M,g)$ tell us about the geometry of $M$.
Consider a pair $(\alpha,n)$ where $\alpha\in \Delta$ and $n\in \ZZ$.
This defines an affine hyperplane
$$   \{H\in \mathfrak{a}\,|\,  \alpha(H) = n \}\subseteq \mathfrak{a}     $$
which we call a \textit{singular plane} and which we denote by $(\alpha,n)$.
As we have seen above, the geodesics in $T$ starting at $o$ are images of lines through the origin in $\mathfrak{a}$ under $\mathrm{Exp}$.
Let $$\sigma_H\colon [0,\infty)\to \mathfrak{a} ,\quad \sigma_H(t) = tH$$
for some $H\in\mathfrak{a}$ be such a ray and let
$$   \gamma_H\colon [0,\infty)\to M,\quad \gamma_H = \Exp\circ \, \sigma_H   $$
be the corresponding geodesic in $M$.
Define times $0 = t_0 < t_1<t_2 < \ldots $ inductively by setting $t_0 = 0$ and
$$    t_i \coloneqq \inf\{ t\in (t_{i-1},\infty)\,|\,   \sigma(t)\text{ lies in a hyperplane }(\alpha,n),\,\,n\neq 0  \}     $$
for $i\geq 1$. Of course, it is also possible that $\sigma$ meets an intersection of several hyperplanes at a time $t_i$, $i\in \NN$.
The following proposition states that the singular planes encode the conjugate points in the symmetric space $M$.
\begin{prop}[\cite{helgason:78}, Proposition VII.3.1] \label{prop_conjugate_points}
Let $\gamma_H$ be defined as above and consider the set $C=\{t_i\}_{i\in \NN}$.
The point $\gamma_H(t)$ is conjugate to $o$ along $\gamma_H$ if and only if $t\in C$.
\end{prop}
We refer to Figure \ref{FigureRootSystem} for a sketch of the singular planes of a symmetric space of rank $2$.
A priori Proposition \ref{prop_conjugate_points} is only a statement about the geodesics in $M$ that start at $p$ and that lie in the maximal torus $T$.
However, the following lemma shows that every geodesic $\gamma\colon[0,\infty)\to M$ is of the form $\gamma_H$ as above up to a global isometry.
Note that the action of $G$ on the symmetric space $M$ induces an action $\Phi\colon G\times \Lambda M\to \Lambda M$ on the free loop space given by
$$ \Phi(g,\gamma)(t) \coloneqq  g\cdot \gamma(t)\quad \text{for}\,\,g\in G,\gamma\in\Lambda M\, .       $$
\begin{lemma} \label{prop_geodesics_symm}
Let $\gamma\colon [0,\infty)\to M$ be a geodesic in $M$.
Then there is an isometry $g\in G$ and an element $H\in\mathfrak{a}$ such that $g \cdot \gamma = \gamma_H$.
\end{lemma}
\begin{proof}
    Since the isometry group acts transitively on $M$, we find an element $g'\in G$ such that $\gamma_1 = g'\cdot \gamma\colon [0,\infty)\to M$ starts at the base point $o = \pi(e)\in M$.
    Let $H_1 = \Dot{\gamma}_1(0)\in \mathfrak{m}$ be the tangent vector at $t= 0$.
    We have 
    $$\gamma_1(t) = \Exp(t\cdot H_1) = \pi\circ \exp(t \cdot H_1) . $$
    However, $H_1$ might not lie in $\mathfrak{a}$.
    By \cite[Theorem VII.8.6]{helgason:78}, we can find an element $k\in K$ such that $H := \mathrm{Ad}_k(H_1) \in \mathfrak{a}$.
    Consider the geodesic $\gamma_2\colon [0,\infty)\to M, t\mapsto k\cdot \gamma_1(t)$.
    We have
    $$   \gamma_2(t) = k\cdot \gamma_1(t) = \pi \big( k \exp(t\cdot H_1) \big) = \pi\big( k \exp(t\cdot H_1 ) k^{-1} \big)  =  \Exp(t\cdot \mathrm{Ad}_k(H_1)) = \gamma_H(t)    $$
    for all $i\in [0,\infty)$.
    Hence, $\gamma_2 = \gamma_H$ and by construction we have $\gamma_2 = (kg')\cdot \gamma$.
    Setting $g = kg'$ yields the claim.    
\end{proof}

Since $\Exp\colon \mathfrak{a}\to \Exp(\mathfrak{a})$ is the universal covering of the torus $T= \Exp(\mathfrak{a})$, we obtain a lattice $\mathcal{F}\subseteq \mathfrak{a}$ by setting
$$   \mathcal{F} = \{ H\in\mathfrak{a}\,|\, \Exp(H) = o\}  \, .     $$
Since the maximal torus $T$ is isometric to a flat torus, the closed geodesics in $T$ are in one-to-one correspondence to points in the lattice $\mathcal{F}$.
See Figure \ref{FigureRootSystem} for a sketch of the lattice of a symmetric space of rank $2$.

By Lemma \ref{prop_geodesics_symm} we see that the critical submanifolds of the energy functional in $\Lambda M$ are the orbits of the closed geodesics in $T$ under the action of $G$.
However, note that several distinct points in the lattice $\mathcal{F}$ can correspond to closed geodesics which lie in the same orbit under the action of $G$.
This indeterminacy can be removed as follows.

Recall from above that we have introduced the set $\Delta$ of roots of $(M,g)$.
One can define certain orderings on $\Delta$, and define the set of \textit{positive roots} $\Delta_+\subseteq \Delta$ with respect to an ordering, see \cite[Section VII.2]{helgason:78} for details. 
We fix such an ordering on $\Delta$.
The \textit{positive Weyl chamber} with respect to this ordering is defined as
$$W:= \{X \in \mathfrak{a} \, | \,  \alpha(X) >0  \text{ for all }  \alpha \in \Delta_+ \}. $$
The positive Weyl chamber is a convex subset of $\mathfrak{a}$. 
For other choices of orderings of $\Delta$ we obtain other convex subsets of $\mathfrak{a}$, which are also called Weyl chambers, but they are disjoint from the positive Weyl chamber.
The \textit{Weyl group} $W(G,K)\subseteq \mathrm{Gl}(\mathfrak{a})$ of $M = G/K$ is the group generated by the reflections $s_{\alpha}\colon \mathfrak{a}\to \mathfrak{a}$ about the hyperplanes 
$$  (\alpha,0) =  \{ H\in \mathfrak{a} \,|\, \alpha(H) = 0\} \,   $$
for all $\alpha\in \Delta$.
The Weyl group is finite and it acts simply transitively on the set of Weyl chambers, see \cite[Proposition VII.2.12]{helgason:78}.
In order to get a one-to-one correspondence between points in the lattice $\mathcal{F}$ and the critical submanifolds we need to consider the intersection of the lattice $\mathcal{F}$ with the closure of the Weyl chamber $\overline{W}$.

\begin{figure}[t]
\centering
\includegraphics[scale=0.55]{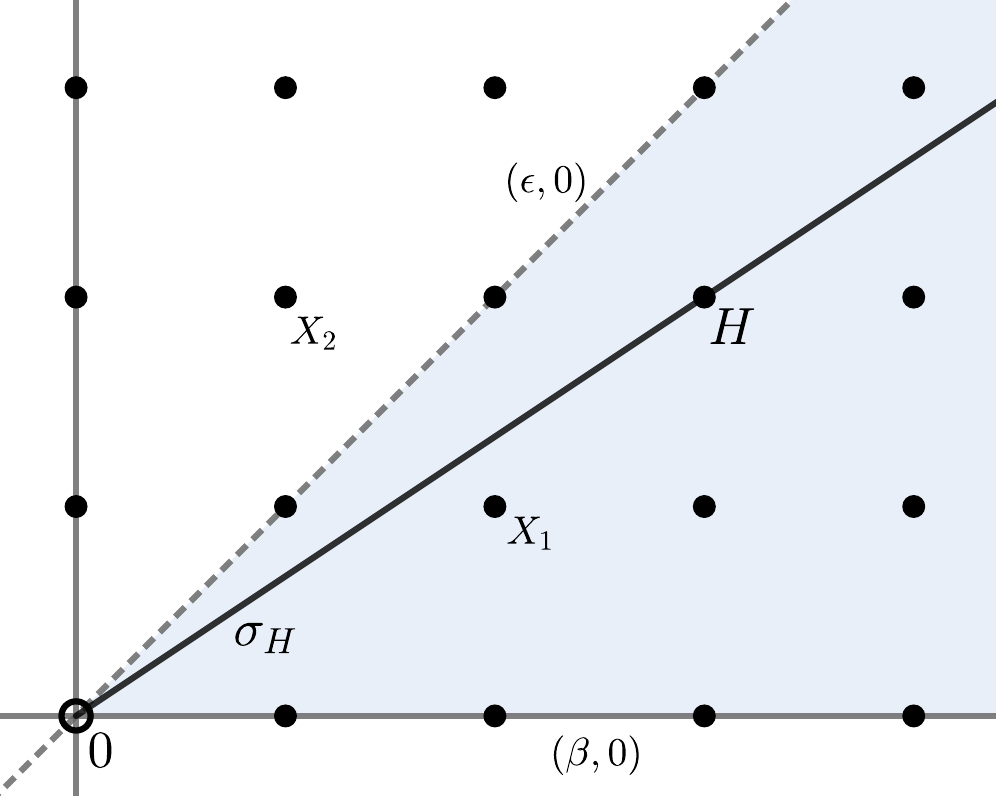}
\caption{Example of the positive Weyl chamber for the case of the maximal abelian subspace shown in Figure \ref{FigureRootSystem}.
The Weyl chamber is depicted as the shaded area.
The lattice points are denoted as bullet points.
Note that the closed geodesic $\gamma_H = \Exp\circ \sigma_H$ is prime, since it intersects no lattice points before reaching $H$.
Further note that the closed geodesics $\gamma_{X_1}$ and $\gamma_{X_2}$ lie in the same critical manifold in $\Lambda M$ since $X_2$ is mapped to $X_1$ by the reflection about the hyperplane $(\epsilon,0)$.}
\label{Figure_weylchamber_lattice}
\end{figure}

\begin{prop}\label{prop_lattice_weyl}
    Let $\gamma\in \Lambda M$ be a closed geodesic.
    Then there is a unique element $H\in \mathcal{F}\cap \overline{W}$ such that $g\cdot \gamma = \gamma_H$ for some $g\in G$.
    Consequently, the critical manifolds of the energy functional in $\Lambda M$ are in one-to-one correspondence with $\mathcal{F}\cap \overline{W}$.
\end{prop}
\begin{proof}
If $\gamma\in \Lambda M$ is a closed geodesic, then by Lemma \ref{prop_geodesics_symm} there is an isometry $g\in G$ of $M$ such that $g\cdot \gamma = \gamma_H$ for some $H\in \mathfrak{a}$.
Moreover, it holds that $H\in\mathcal{F}$, since $\gamma$ is a closed geodesic and $g\cdot \gamma$ is thus a closed geodesic as well.
By \cite[Theorem VII.2.22]{helgason:78} there exists a unique element $H'\in \overline{W}$ which is in the orbit of $H$ under the action of the Weyl group.
Since the action of the Weyl group is lattice-preserving, see \cite[Proof of Theorem 8.5]{helgason:78}, we have $H'\in \mathcal{F}\cap \overline{W}$.
Hence, there is an element $s\in W(G,K)$ of the Weyl group such that $s (H) = H'$.
By \cite[Corollary VII.2.13]{helgason:78} there is a $k_0\in K$ such that $s = \mathrm{Ad}_{k_0}|_{\mathfrak{a}}$.
Using the same computation as in the proof of Lemma \ref{prop_geodesics_symm} we see that $(k_0\cdot g)\cdot \gamma = \gamma_{H'}$.
\end{proof}
For a sketch of the positive Weyl chamber and the lattice, see Figure \ref{Figure_weylchamber_lattice}.
Recall that we call a closed geodesic $\gamma\in \Lambda M$ \textit{prime} if it is not the iterate of another closed geodesic, i.e. if there is no $l\geq 2$ and a closed geodesic $\sigma\in \Lambda M$ with $\gamma = \sigma^l$.
We say that a critical submanifold is prime if it consists only of prime closed geodesics.
The next proposition says that there are many prime critical submanifolds for higher rank symmetric spaces.
\begin{prop}\label{infty_many_prime_geodesics}
Let $M$ be a compact symmetric space. If the rank of $M$ satisfies $r\geq 2$, there are countably infinite distinct prime critical submanifolds in $\Lambda M$.
\end{prop}
\begin{proof}
The lattice $\mathcal{F}\subseteq \mathfrak{a}$ is generated by elements $X_1,\ldots, X_r$ where $r$ is the rank of $M$.
Every element of the form
\begin{equation} \label{lattice_prime}
       \lambda_1 X_1 + \ldots + \lambda_r X_r  \quad \text{with at least one } \lambda_i = \pm 1  , \,\, i\in\{1,\ldots,r\}  
\end{equation}
corresponds to a prime closed geodesic.
If $r\geq 2$ there are countably infinite such elements.
Recall that all Weyl chambers are homeomorphic to one another and that the union of the closures of the finitely many Weyl chambers is all of $\mathfrak{a}$.
Since the action of the Weyl group preserves the lattice $\mathcal{F}$, see \cite[Theorem VII.8.5]{helgason:78}, there must be countably infinite elements of the form \eqref{lattice_prime} in each Weyl chamber.
\end{proof}
In contrast, for a rank one symmetric space it is well-known that there is exactly one prime critical submanifold of $\Lambda M$.

\medskip
We conclude this section by computing the index growth of the closed geodesics in the compact symmetric space $M$.
We have seen above that conjugate points appear along a closed geodesic of the form $\gamma_H = \Exp\circ \,\sigma_H$ for $H\in\mathfrak{a}$ whenever the ray $\sigma_H$ intersects a singular plane.
We shall now study an alternative characterization of the conjugate points along closed geodesics which will be used in the next section.
Let $c$ be a closed geodesic in $M$ starting at the basepoint $o\in M$.
Let $K_c$ be the closed subgroup of $K$ that keeps $c$ pointwise fixed, i.e.
$$   K_c  =  \{ k\in K\,|\,   k . c(t) = c(t) \,\,\text{for all}\,\, t\in I\}\,  .        $$
The conjugate points along $c$ can be characterized as follows.
A point $q = c(t) \in M$ with $t\in I$ is a conjugate point along $c$ if and only if the stabilizer 
$$      K_q \coloneqq \{ k\in K\, |\, k . q = q\}       $$
of $q$ satisfies 
$$     \mathrm{dim}(K_q) > \mathrm{dim}(K_c)\,  .    $$
This can be seen as follows.
Let $\mathfrak{k}_q$ be the Lie algebra of $K_q$ and $\mathfrak{k}_c$ be the Lie algebra of $K_c$.
The Lie algebra $\mathfrak{k}_c$ is a Lie subalgebra of $\mathfrak{k}_q$.
Let $\xi\in\mathfrak{k}_q$ be an element such that $\xi\not\in \mathfrak{k}_c$ and consider the one-parameter group $\Xi\colon \mathbb{R}\to K_q, \Xi(s) = \exp(s\xi)$.
The variation of the geodesic $c$ defined by $$[0,t]\times (-\epsilon,\epsilon)\to M , \quad (\tau,s)\mapsto  \exp(s\xi) \cdot c(\tau)  $$
is a non-trivial variation of $c$ through geodesics with fixed endpoints.
It thus determines a non-trivial Jacobi field and hence $q = c(t)$ is conjugate to $o$ along $c$.
Elaborating on this idea one finds that the difference $\mathrm{dim}(K_q) - \mathrm{dim}(K_c)$ is precisely the multiplicity of the conjugate point $q$, i.e. the dimension of the space of Jacobi fields along $c|_{[0,t]}$ that vanish at the endpoints, see \cite[p. 13-14]{ziller:1977}.

Recall that the closed geodesics in $M$ are the critical points of the energy functional on $\Lambda M$.
We want to express the Morse index of a closed geodesic $c\in \Lambda M$ through the geometric data at hand.
One shows that in the case of $(M,g)$ being a symmetric space, the index of the closed geodesic $c$ is equal to the sum of the multiplicities of the conjugate points in the interior, see \cite[Theorem 4]{ziller:1977}.
This means that we have
$$      \mathrm{ind}(c) = \sum_{0<t_i<1} (\mathrm{dim}(K_{c(t_i)}) - \mathrm{dim}(K_c))\, ,$$
where the times $t_i\in(0,1)$ are such that $c(t_i)$ is a conjugate point as above.
Define the integer 
$$\mu =\mathrm{dim}(K) -\mathrm{dim}(K_c)\,.
$$
If $c$ is the $k$-th iterate of a prime closed geodesic $\sigma$, then by a simple counting argument we obtain
\begin{equation} \label{eq_index}
      \mathrm{ind}(c) = k \, \mathrm{ind}(\sigma) +  (k-1) \mu  
\end{equation}
since the origin appears $(k-1)$ times in the interior of $c$ and clearly has stabilizer $K$.

We now consider the critical manifolds of the energy functional in $\Lambda M$.
By Proposition \ref{prop_lattice_weyl} the critical manifolds are in one-to-one correspondence with the lattice points that intersect the closure of the positive Weyl chamber, i.e. with the set $\mathcal{F}\cap \overline{W}$.
Note that for the connected component $\Sigma_c$ of the critical set of $E$ which contains the closed geodesic $c= \gamma_H$ for $H\in\mathcal{F}\cap \overline{W}$, the $G$-action on $\Lambda M$ induces a transitive $G$-action on $\Sigma_c$.
This is again due to Proposition \ref{prop_lattice_weyl}.
The stabilizer group of $c$ under this action is the group $K_c$ that we introduced above.
Consequently, the critical manifold $\Sigma_c$ is diffeomorphic to the homogeneous space $G/K_c$.
Ziller has shown that $G/K_c$ is a non-degenerate critical submanifold of $E$, see \cite[Theorem 2]{ziller:1977}.
Since this holds for each critical submanifold, $E$ is a Morse-Bott function.
The nullity of the closed geodesic $c$, i.e. the dimension of the kernel of the Hessian of $E$ at $c$, is thus equal to the dimension of the homogeneous space $G/K_c$.
Since the canonical map 
$$G/K_c \to G/K,\quad g K_c \mapsto gK $$ 
is a fiber bundle with typical fiber $K/K_c$, we see that
$$  \mathrm{null}(c) = \mathrm{dim}(G/K) + \mathrm{dim}(K/K_c)  = n + \mu \, .     $$
Consequently, the sum of index and nullity of the $k$-th iterate of a prime closed geodesic $\sigma$ is
\begin{equation} \label{eq_index_null}
       (\mathrm{ind} + \mathrm{null})(c) = k\,\mathrm{ind}(\sigma) + k \,\mu + n\,  . 
\end{equation}

Note that the homogeneous space $K/K_c$ is diffeomorphic to the orbit $K . c\subseteq \Omega_o M$ of the closed geodesic $c$ under the induced action of $K$ on the based loop space.
Moreover, one can also identify $K_c$ with the subgroup
$$   \{ k\in K\,|\, \Ad_k(\Dot{c}(0)) = \Dot{c}(0)\} \subseteq K  \,   , $$
see \cite[p. 14]{ziller:1977}.
The group $K$ acts by linear isometries on the space $\mathfrak{m}$, hence the orbit of $\Dot{c}(0)$ is diffeomorphic to the space $K/K_c$ which is embedded in the hypersphere
$$   S_{c} \coloneqq  \{ X\in \mathfrak{m}\,|\, |X| = |\Dot{c}(0)|\} \cong \mathbb{S}^{n-1}\,  .      $$
Recall that the dimension of $K/K_c$ is $\mu$.
As usual we can assume that $\Dot{c}(0)\in \mathfrak{a}$.

If $M$ is an $n$-dimensional symmetric space of rank $1$, then the space $K/K_c$ is the whole sphere $S_c$, since the compact rank one symmetric spaces are precisely the symmetric spaces where all geodesics are closed and of common prime length.
Thus $\mu = n-1$ in this case.
If $M$ is $n$-dimensional and of rank $r\geq 2$, then the orbit of $K$ cannot be the whole hypersphere $S_c\subseteq \mathfrak{m}$, which one can see as follows.
If the orbit of $K$ was the whole $S_c$ then in particular every point in $S_c\cap \mathfrak{a}$ lies in the orbit $K . \Dot{c}(0)$.
But if $H\in\mathfrak{a}$ is in the $K$-orbit of $\Dot{c}(0)$, i.e. there is a $k\in K$ with $H = \Ad_k(\Dot{c}(0))$, then there is an element of the Weyl group $s\in W(G,K)$ with $H = s (\Dot{c}(0))$, see \cite[Proposition VII.2.2]{helgason:78}.
But there are only finitely many elements of the Weyl group, so for rank $r\geq 2$ this yields a contradiction to the assumption that the $K$-orbit of $\Dot{c}(0)$ is all of $S_c$.
Consequently, $K/K_c$ is an embedded closed submanifold of $\mathbb{S}^{n-1}$, which is not the whole sphere and therefore must have positive codimension.
We have thus shown the following.
\begin{lemma}\label{remark_mu}
    Let $M$ be a symmetric space of rank $r$ and $c\in \Lambda M$ be a closed geodesic.
    The dimension $\mu = \mathrm{dim}(K/K_c)$ of the orbit of $c$ under the action of $K$ is equal to $n-1$ if and only if $r = 1$ and $\mu< n-1$ if $r\geq 2$.
\end{lemma}
Recall from above that the number $\mu$ shows up in the formulae for index and index plus nullity of the iterates of $c$.
In Remark \ref{remark_index-growth} we will relate the growth of the index of iteration of closed geodesics in compact symmetric spaces to the results of the next sections concerning the string topology operations on compact symmetric spaces.

\section{The Chas-Sullivan product on symmetric spaces} \label{sec_chas_sullivan}

Let $M= G/K$ be a compact symmetric space.
As we have seen in the previous section, $M$ is a homogeneous space of the connected component of its isometry group $G$ with $K$ being the stabilizer subgroup of a point. In this section we are going to show that each critical submanifold $\Sigma\subseteq \Lambda M$ consisting of prime closed geodesics of $M$ gives rise to a non-nilpotent element in the Chas-Sullivan algebra of $\Lambda M$. 
Throughout the latter part of this section we only consider homology with $\mathbb{Z}_2$-coefficients.
At the end of this section we shall discuss coefficients in an arbitrary unital commutative ring.

Let $\gamma$ be a closed geodesic in $M$.
Let $o = \pi(e)$ be the image of the neutral element $e\in G$ under the canonical projection $\pi \colon G\to M = G/K$.
We choose $o$ as the basepoint of $M$.
Recall from Section \ref{sec_symm_geodesics} that there is an element $g\in G$ such that $c = g.\gamma \in \Lambda M$ satisfies $c(0) = o$.
Since $G$ acts by isometries on $M$, the loop $c$ is also a closed geodesic. 
Therefore the critical manifold of $c$ in $\Lambda M$ with respect to the energy functional is diffeomorphic to $\Sigma_c = G/K_c$ where $K_c$ is the isotropy group
$$   K_c \coloneqq \{ k\in K\,|\, k.c(t) = c(t)\,\,\text{for all}\,\,t\in [0,1]\} \subseteq K\,.    $$
The projection $\Sigma_c = G/K_c \to G/K$ is a fibre bundle projection 
$$K/K_c \hookrightarrow \Sigma_c \xrightarrow[]{\pi} M$$
with homogeneous fibre $K/K_c $.

We shall now describe Ziller's explicit cycles, defined as completing manifolds for the critical manifolds $\Sigma_c$ in \cite[Section 3]{ziller:1977}.
Fix a closed geodesic $c\in \Lambda M$ with $c(0) = o$. Again from Section \ref{sec_symm_geodesics} we recall the following: If $q\in M$ is a point we denote its stabilizer under the action of $K$ by $K_q$, i.e.
$$  K_q = \{k\in K\,|\, k.q = q\}\, .    $$
As Ziller \cite{ziller:1977} and Bott and Samelson \cite{bott:1958a} argue, there are finitely many times $0<t_1<\ldots <t_l< 1$ such that $\mathrm{dim}(K_{c(t_i)}) > \mathrm{dim}(K_c)$ for $i\in\{1,\ldots,l\}$.
The points $c(t_i)$, $i\in\{1,\ldots,l\}$ are precisely the conjugate points along $c$ and the multiplicity of such a conjugate point is equal to $\mathrm{dim}(K_{c(t_i)}) - \mathrm{dim}(K_c)$.
Moreover, the Morse index of the closed geodesic $c$ can be expressed as
$$   \mathrm{ind}(c) = \sum_{i= 1}^l \big(  \mathrm{dim}(K_{c(t_i)}) - \mathrm{dim}(K_c)   \big) \, .  $$
If $c$ has multiplicity $m$, i.e. there is an underlying prime closed geodesic $\sigma$ with $c=\sigma^m$, then it holds that $c(\tfrac{j}{m}) = o$ for $j\in\{1,\ldots,m-1\}$ and the corresponding stabilizer is the whole group $K$.
As before, define $\mu \coloneqq\mathrm{dim}(K) - \mathrm{dim}(K_c)$.
Recall from Section \ref{sec_symm_geodesics} that the index of the iterates of $\sigma$ equals
$$ \mathrm{ind}(c)=  \mathrm{ind}(\sigma^m) = m \mathrm{ind}(\sigma) + (m-1)\mu \,.   $$
Moreover, we have
$$  \mathrm{dim}(\Sigma_c) = n + \mu\,,   $$
where $n$ is the dimension of $M$.
We consider the product of groups
$$  W_c \coloneqq G\times K_{c(t_1)}\times K_{c(t_2)} \times \ldots \times K_{c(t_l)}\, .     $$
There is a right-action of the $(l+1)$-fold product of $K_c$ on $W_c$ given by 
\begin{eqnarray} \label{eq_definition_chi}   \chi : W_c \times K_c^{l+1} &\to& W_c      \\
                ( (g_0,k_1,\ldots,k_l), (h_0,\ldots,h_l)) &\mapsto& (g_0 h_0,h_0^{-1} k_1 h_1,\ldots, h_{l-1}^{-1} k_l h_l)\,  . \nonumber
\end{eqnarray}
This is a proper and free action and we denote the quotient by
$$  \Gamma_c \coloneqq W_c/ (K_c^{l+1})  \, . $$
Note that there is a submersion
\begin{equation}\label{gamma_to_sigma}
  p\colon \Gamma_c  \to \Sigma_c\cong G/K_c \,, \quad p([(g_0,k_1,\ldots,k_l)]) = [g_0]\, .           
\end{equation}

Furthermore there is an embedding $\Sigma_c\to \Gamma_c$ defined by
$$   s\colon \Sigma_c\to\Gamma_c\,, \quad s([g]) = [(g,e,\ldots,e)]     $$
and one sees that $p\circ s = \mathrm{id}_{\Sigma_c}$.

We now describe an embedding $ \Gamma_c \hookrightarrow \Lambda M$.
Define $\widetilde{f}_c\colon W_c \to \Lambda M$ by
$$   \widetilde{f}_c(g_0,k_1,\ldots,k_{l})(t)  = \begin{cases} g_0. c(t), & 0\leq t\leq t_1 
\\
g_0k_1.c(t), & t_1 \leq t \leq t_2
\\
\,\,\,\,\, \vdots & \,\,\,\,\, \vdots
\\
g_0k_1\ldots k_{l}. c(t), & t_l \leq t \leq 1 \,\,\,.
\end{cases}       $$ 
This factors through the action of $K_c^{l+1}$ and therefore defines a map $f_c\colon \Gamma_c\to \Lambda M$.
It can be seen directly that this is an embedding and that $f_c\circ s\colon \Sigma_c\to \Lambda M$ embeds $\Sigma_c$ precisely as the critical manifold.
Moreover, the only closed geodesics in $f_c(\Gamma_c)$ are in this critical manifold.
Hence, we see that $\Gamma_c$ is a completing manifold if coefficients are taken in $\ZZ_2$. 
The same is true with $R$-coefficients if $\Sigma_c$ and $\Gamma_c$ are $R$-orientable as follows from Section \ref{sec_completing}.

We now show that Ziller's cycles behave well with respect to the concatenation map 
$$  \mathrm{concat}_s \colon \Lambda\times_M\Lambda \to \Lambda\,, \quad  (\gamma_1,\gamma_2)\mapsto \gamma\,,     $$
where $s\in(0,1)$ and where 
$$   \gamma(t) \coloneqq \begin{cases} 
\gamma_1(\tfrac{t}{s} ) & 0\leq t\leq s \\ \gamma_2(\tfrac{t+s}{1-s}) & s\leq t\leq 1 \,\,\,.
\end{cases}
$$
Let $\sigma\in \Lambda M$ be a prime closed geodesic and suppose that $\gamma_1,\gamma_2\in \Lambda M$ are iterates of $\sigma$ with multiplicities $m_1$ and $m_2$.
Just as above one can set up completing manifolds $\Gamma_1$ and $\Gamma_2$ corresponding to the critical manifolds $\Sigma_1$ and $\Sigma_2$ where $\Sigma_i \subseteq \Lambda M$ is the critical manifold containing $\gamma_i$ for $i\in\{1,2\}$.
We want to study the relation between $\Gamma_1,\Gamma_2$ and the completing manifold $\Gamma_3$ associated to the concatenation $\gamma_3 = \sigma^{m_1+m_2}$.
We set $m_3 = m_1 + m_2$ and note that we have
$$  \Sigma_3 \cong \Sigma_1 \cong \Sigma_2 \cong G/K_{\sigma}\, .     $$
Recall that we have submersions $p_i\colon \Gamma_i\to \Sigma_i$, $i\in\{1,2,3\}$, and that there is a fiber bundle $\pi\colon G/K_{\sigma}\to M$.
This gives submersions $\pi\circ p_i\colon \Gamma_i\to M$ and we can look at the fiber product
$$  \Gamma_1\times_M \Gamma_2   \subseteq \Gamma_1\times \Gamma_2  $$
with respect to $\pi\circ p_1$ and $\pi\circ p_2$.
Clearly, $\Gamma_1\times_M\Gamma_2$ is a submanifold of $\Gamma_1\times\Gamma_2$ of codimension $n$ and it has a canonical map
$$  p_1\times p_2\colon \Gamma_1\times_M \Gamma_2 \to \Sigma_1\times_M \Sigma_2,   \quad (p_1\times p_2)(X,Y) = (p_1(X),p_2(Y))\, .   $$
Moreover, we see that 
\begin{equation} \label{fiber_product_gamma}
      \Gamma_1\times_M\Gamma_2 = \{ ([(g_0,k_1,\ldots,k_{l_1})],[(\widetilde{g}_0,\widetilde{k}_1,\ldots,\widetilde{k}_{l_2})])\in \Gamma_1\times\Gamma_2\,|\, g_0 K = \widetilde{g}_0 K\} \,. 
\end{equation}
Recall the construction of the group $W_c$ above for a closed geodesic $c$.
If we define $H$ by $ W_{\sigma} = G\times H$ where $\sigma$ is the underlying prime closed geodesic of $\gamma_1$ and $\gamma_2$ as before, then we see that
$$  W_{\gamma_i} = G\times H \times (K\times H)^{m_i-1}\,,\quad i\in\{1,2,3\}\, .     $$

Define a map $\Phi\colon \Gamma_1\times_M \Gamma_2\to \Gamma_3$ by
$$  \Phi([( g_0,k_1,\ldots,k_{l_1})],[(\widetilde{g}_0,\widetilde{k}_1,\ldots,\widetilde{k}_{l_2})]) = [(g_0,k_1,\ldots,k_{l_1} ,  (g_0 k_1\ldots k_{l_1})^{-1} \widetilde{g}_0,  \widetilde{k}_1,\ldots,\widetilde{k}_{l_2})]   . $$
Note that the map $\Phi$ is smooth and indeed well-defined since $g_0^{-1}\widetilde{g}_0\in K$ by equation \eqref{fiber_product_gamma}.
Set $\tau = \tfrac{m_1}{m_1+m_2}$ and define $\phi\colon \Sigma_1\times_M\Sigma_2\to \Sigma_3$ by
$$   \phi([g],[\widetilde{g}]) = [g] , $$
where we identify $\Sigma_i\cong G/K_{\sigma}$ as usual.

\begin{lemma} \label{completing_manifolds_diffeo}
The map $\Phi$ is a diffeomorphism.
Moreover, it is compatible with the completing manifold structure in the sense that the diagrams
$$
\begin{tikzcd}
\Gamma_1\times_M \Gamma_2 \arrow[]{r}{f_1\times f_2} 
\arrow[d , "\Phi", "\cong"']  
& [2em]  \Lambda\times_M \Lambda \arrow[d, "\mathrm{concat}_{\tau}" ]
\\
\Gamma_3 \arrow[]{r}{f_3} & \Lambda
\end{tikzcd}
$$
and 
$$
\begin{tikzcd}
\Gamma_1\times_M \Gamma_2 \arrow[]{r}{p_1\times p_2} 
\arrow[d , "\Phi", "\cong"']  
&  [2em]  \Sigma_1\times_M \Sigma_2 \arrow[d, "\phi" ]
\\
\Gamma_3 \arrow[]{r}{p_3} & \Sigma_3
\end{tikzcd}
$$
commute.
\end{lemma}
\begin{proof}
An explicit inverse for $\Phi$ is given by $\Psi\colon \Gamma_3\to \Gamma_1\times_M\Gamma_2$,
$$  \Psi([(g_0,k_1,\ldots,k_{l_3})]) = ([(g_0,k_1,\ldots,k_{l_1})],[((g_0 k_1 \ldots k_{l_1}) k_{l_1 +1} , k_{l_1 + 2},\ldots,k_{l_3}) ]) \,  .   $$
Again, one checks that this is well-defined and smooth and is in fact an inverse to $\Phi$.
The commutativity of the diagrams can then be checked by unwinding the definitions.
\end{proof}

The property of the completing manifolds described in the following corollary is analogous to the property of the completing manifolds constructed by Oancea for spheres and complex projective spaces in \cite[Section 6]{oancea:2015}.

\begin{cor}
Let $\sigma$ be a prime closed geodesic of $M$. Let $\Gamma$ denote the completing manifold of the prime closed geodesics isometric to $\sigma$ and let $\Gamma^m$ denote the completing manifold corresponding to $\sigma^m$, for $m\in\NN$. We have
\begin{equation*}
    \Gamma^m\cong\underbrace{\Gamma\times_M\Gamma\times_M\dots\times_M\Gamma}_{m \text{ times}}\,.
\end{equation*}
\end{cor}

In the following we use the notation $(Y,\sim A) \coloneqq (Y, Y\setminus A)$ for a pair of topological spaces $(Y,A)$ and $\Lambda\coloneqq\Lambda M$.
Recall from Section \ref{sec_def} that for two homology classes $X\in \mathrm{H}_i(\Lambda M)$ and $Y\in \mathrm{H}_j(\Lambda M)$ their Chas-Sullivan product $X\wedge Y$ is defined by the composition
\begin{equation*}
    \xymatrix{
    \mathrm{H}_i(\Lambda)\otimes \mathrm{H}_j(\Lambda )\ar[d]_-{(-1)^{n(n-j)}\times}&& \mathrm{H}_{i+j-n}(\Lambda )\\
    \mathrm{H}_{i+j}(\Lambda \times\Lambda )\ar[r]&\mathrm{H}_{i+j}(U_{CS},\sim \Lambda \times_M\Lambda )\ar[r]^-\cong&\mathrm{H}_{i+j-n}(\Lambda \times_M\Lambda)\ar[u]_-{\mathrm{concat}_*}
    }
\end{equation*}
where $U_{CS}$ is a tubular neighbourhood for the figure eight space $\Lambda \times_M\Lambda\subset\Lambda \times\Lambda$ and the isomorphism in the bottom line is a Thom isomorphism. Consider the following commutative diagram
\begin{equation*}
    \xymatrix{
    &\Lambda\times_M \Lambda \ar@{^{(}->}[rr]
&&\Lambda\times \Lambda \ar@/^2pc/[lddd]^-{\mathrm{ev}_0\times \mathrm{ev}_0}\\
 \Gamma_1\times_M \Gamma_2\ar[d]\ar@{^{(}->}[rr]\ar@{^{(}->}[ru]\ar[d]&&\Gamma_1\times \Gamma_2\ar[d]_-{p_1\times p_2}\ar@{^{(}->}[ru]^-{f_1\times f_2}& \\
 \Sigma_1\times_M \Sigma_2\ar[d]\ar@{^{(}->}[rr]&&\Sigma_1\times \Sigma_2\ar[d]_-{\pi\times\pi}&\\
 M\ar@{^{(}->}[rr]^-\Delta&&M\times M&
    }
\end{equation*}
where the arrows on the left hand side are given by restricting the ones on the right hand side of the diagram.
Since $\mathrm{ev}_0\times \mathrm{ev}_0$, $p_1\times p_2$ and $\pi\times\pi$ are submersions it follows by transversality that
\begin{align*}
   (\mathrm{ev}_0\times \mathrm{ev}_0)^{-1}\left(\Delta(M)\right)&=\Lambda\times_M \Lambda\hookrightarrow\Lambda\times \Lambda\,\\
   (f_1\times f_2)^{-1}(\Lambda\times_M \Lambda)=\left((\pi\times\pi)\circ(p_1\times p_2)\right)^{-1}\left(\Delta(M)\right)&=\Gamma_1\times_M \Gamma_2\hookrightarrow\Gamma_1\times \Gamma_2\,\\
   (\pi\times \pi)^{-1}\left(\Delta(M)\right)&=\Sigma_1\times_M \Sigma_2\hookrightarrow\Sigma_1\times \Sigma_2\,.
\end{align*}

are inclusions of smooth submanifolds of codimension $n$. Let $N_{\Sigma_1\times_M \Sigma_2}$, $N_{\Gamma_1\times_M \Gamma_2}$, $N_{\Lambda\times_M \Lambda}$ denote the corresponding normal bundles. If $N_M$ denotes the normal bundle of the diagonal inclusion $\Delta\colon M\to M\times M$ it also follows that
\begin{align*}
    (\mathrm{ev}_0\times \mathrm{ev}_0)^*(N_M)&\cong N_{\Lambda\times_M \Lambda}\,\\
    (f_1\times f_2)^*(N_{\Lambda\times_M \Lambda})\cong(f_1\times f_2)^*\circ(\mathrm{ev}_0\times \mathrm{ev}_0)^*(N_M)\cong (p_1\times p_2)^*\circ (\pi\times\pi)^*(N_M)&\cong N_{\Gamma_1\times_M \Gamma_2}\,\\
    (\pi\times\pi)^*(N_M)&\cong N_{\Sigma_1\times_M \Sigma_2}\,.
\end{align*}
holds. Furthermore, it holds that if $U_M$ is a tubular neighbourhood of $\Delta\colon M\to M\times M$, then 
\begin{align*}
(\mathrm{ev}_0\times \mathrm{ev}_0)^{-1}\left(U_M\right)=U_{CS}&\text{ is a tubular neighbourhood of }\Lambda\times_M \Lambda\hookrightarrow\Lambda\times \Lambda\,,\\
(f_1\times f_2)^{-1}(U_{CS})\eqqcolon U&\text{ is a tubular neighbourhood of }\Gamma_1\times_M \Gamma_2\hookrightarrow\Gamma_1\times \Gamma_2\,.
\end{align*}
It then also follows that the diagram
\begin{equation*}
\xymatrixcolsep{2.5pc}
    \xymatrix{
    \mathrm{H}_{*}(U,\sim \Gamma_1\times_M \Gamma_2)\ar[r]^-{(f_1\times f_2)_*}\ar[d]_-{ \cap (f_1\times f_2)^*\tau}^\cong&\mathrm{H}_{*}(U_{CS},\sim \Lambda \times_M\Lambda )\ar[d]^-{ \cap \tau}_-\cong\\
    \mathrm{H}_{* -n}(U)\ar[r]^-{(f_1\times f_2)_*}\ar[d]^-\cong  &   \mathrm{H}_{* -n}(U_{CS})\ar[d]_-\cong\\
    \mathrm{H}_{* -n}(\Gamma_1\times_M \Gamma_2)\ar[r]^-{(f_1\times f_2)_*} &  \mathrm{H}_{* -n}(\Lambda \times_M\Lambda )
    }
\end{equation*}
is commutative, where the vertical compositions are the Thom isomorphisms. This can for instance be seen by using a Riemannian submersion metric on $\Gamma_1\times\Gamma_2\to M\times M$. Using Lemma \ref{completing_manifolds_diffeo} this commutative diagram inserts into the following larger commutative diagram:
\begin{equation*}
\xymatrixcolsep{2.5pc}
    \xymatrix{
\mathrm{H}_i(\Gamma_1)\otimes \mathrm{H}_j(\Gamma_2)\ar[d]^-\times\ar[r]^-{(f_1)_*\otimes (f_2)_*} & \mathrm{H}_i(\Lambda)\otimes \mathrm{H}_j(\Lambda)\ar[d]^-\times\\
\mathrm{H}_{i+j}(\Gamma_1\times\Gamma_2)\ar[d]\ar[r]^-{(f_1\times f_2)_*} & \mathrm{H}_{i+j}(\Lambda\times\Lambda)\ar[d]\\
\mathrm{H}_{i+j}(U,\sim\Gamma_1\times_M \Gamma_2)\ar[d]^\cong\ar[r]^-{(f_1\times f_2)_*} & \mathrm{H}_{i+j}(U_{CS},\sim \Lambda\times_M \Lambda)\ar[d]_\cong\\
\mathrm{H}_{i+j-n}(\Gamma_1\times_M \Gamma_2)\ar[d]_{\Phi_*}^\cong\ar[r]_-{(f_1\times f_2)_*} & \mathrm{H}_{i+j-n}(\Lambda\times_M \Lambda)\ar[d]^-{\mathrm{concat}_*}\\
\mathrm{H}_{i+j-n}(\Gamma_3)\ar[r]_-{(f_3)_*} & \mathrm{H}_{i+j-n}(\Lambda)
    }
\end{equation*}
Here the vertical composition on the right is the Chas-Sullivan product $\wedge$ up to sign.

From now on we only consider homology with $\mathbb{Z}_2$-coefficients.

\begin{theorem}\label{cs_prod_orientation_classes}
The class $(f_1)_*([\Gamma_1])\wedge(f_2)_*([\Gamma_2])$ is nonzero and equal to $(f_3)_*([\Gamma_3])$.
\end{theorem}

\begin{proof}
We know from Section \ref{sec_completing} that the orientation class $[\Sigma_c]$ of a completing manifold of a closed geodesic $c$ injects via $(f_c)_*$ into $\mathrm{H}_*(\Lambda)$. 
More precisely, one sees that $[\Gamma_c]= (p_c)_!\left([\Sigma_c]\right)$, where $p_c\colon \Gamma_c\to\Sigma_c$ is the projection map from above, see equation \eqref{gamma_to_sigma}, and consequently
$$    (f_c)_*[\Gamma_c]  = (f_c)_*\circ (p_c)_! [\Sigma_c]    $$
is non-zero.
We have to show that composition on the left of the above diagram maps $[\Gamma_1]\otimes[\Gamma_2]$ to $[\Gamma_3]$. 
This can be seen as follows. 
By the Künneth theorem we have that $[\Gamma_1]\times[\Gamma_2]$ equals $[\Gamma_1\times\Gamma_2]$. 
The Thom class of $\Gamma_1\times_M\Gamma_2\subset\Gamma_1\times\Gamma_2$ is the Poincaré dual of $[\Gamma_1\times_M\Gamma_2]$ in the homology of $(\bar{U},\partial\bar{U})$ where $\bar{U}$ is the closure of the tubular neighborhood $U$, see \cite[Section VI.11]{bredon:2013}.
Consequently, it follows that capping with the Thom class maps $[\Gamma_1\times\Gamma_2]$ to $[\Gamma_1\times_M\Gamma_2]$.
Since $\Phi$ is a diffeomorphism the orientation class of $\Gamma_1\times_M \Gamma_2$ is mapped to the orientation class of $\Gamma_3$.
\end{proof}

\begin{cor} \label{cor_non_nilpotent}
Let $\sigma$ be a prime closed geodesic, $\Sigma$ the corresponding critical submanifold and $p_\sigma:\Gamma\to\Sigma$ the projection of the completing manifold onto $\Sigma$. Then the class
\begin{equation*}
  \Theta_\sigma\coloneqq (f_\sigma)_*\circ (p_\sigma)_!\left([\Sigma]\right)=(f_\sigma)_*\left([\Gamma]\right)\,.  
\end{equation*}
is non-nilpotent in the Chas-Sullivan algebra of $\Lambda$.
\end{cor}

We now use that $\Sigma_1$ and $\Sigma_2$ are diffeomorphic. Hence, we just write $\Sigma$ to denote either of them. There are even more non-zero Chas-Sullivan products stemming from non-zero products in the intersection algebra of $\Sigma\cong G/K_\sigma$:

\begin{theorem} \label{theorem_cs_intersection}
If two classes $a,b\in \mathrm{H}_*(\Sigma)$ have nonzero intersection product, then
\begin{equation*}
    (f_1)_*\circ(p_1)_!(a)\wedge(f_2)_*\circ(p_2)_!(b)\neq 0\,.
\end{equation*}
\end{theorem}

\begin{proof}
Consider the inclusions
\begin{equation*}
    \xymatrix{
    \Sigma\ar[r]^-{\Delta'}&\Sigma\times_M\Sigma\ar[r]^-J&\Sigma\times\Sigma
    }
\end{equation*}
whose composition gives the diagonal embedding $\Delta_{\Sigma}\colon\Sigma\to\Sigma\times\Sigma$. The map $\Delta'$ is the diagonal section of the map $\phi$ given in Lemma \ref{completing_manifolds_diffeo}. The intersection product $a\bullet b$ of $a,b$ is given by
\begin{align*}
   a\bullet b&= \mathrm{PD}_{\Sigma}\circ(\Delta_{\Sigma})^* \circ\mathrm{PD}_{\Sigma\times\Sigma}^{-1}(a\times b)=   \mathrm{PD}_{\Sigma}\circ(\Delta')^* \circ J^* \circ\mathrm{PD}_{\Sigma\times\Sigma}^{-1}(a\times b)\\
  &=(\Delta')_!\circ J_!(a\times b)\,.
\end{align*}
Let $\lambda_1\coloneqq \mathrm{ind}(\sigma^{m_1})$ and $\lambda_2\coloneqq \mathrm{ind}(\sigma^{m_2})$ denote the Morse indices of the iterates $\gamma_1=\sigma^{m_1}$ and $\gamma_2=\sigma^{m_2}$ of $\sigma$. The intersection product $a\bullet b$ coincides with left most vertical composition in the following extension of the diagram above, compare Appendix B of \cite{hingston:2017}:
\begin{equation*}
\xymatrixcolsep{2.5pc}
    \xymatrix{
    \mathrm{H}_{i-\lambda_1}(\Sigma)\otimes \mathrm{H}_{j-\lambda_2}(\Sigma)\ar[d]_-\times\ar[r]^-{(p_1)_!\otimes (p_2)_!} & \mathrm{H}_i(\Gamma_1)\otimes \mathrm{H}_j(\Gamma_2)\ar[d]^-\times\ar[r]^-{(f_1)_*\times (f_2)_*} & \mathrm{H}_i(\Lambda)\otimes \mathrm{H}_j(\Lambda)\ar[d]^-\times\\
    \mathrm{H}_{i+j-\lambda_1-\lambda_2}(\Sigma\times\Sigma)\ar[r]^-{(p_1\times p_2)_!}\ar[dd]_-{J_!} & \mathrm{H}_{i+j}(\Gamma_1\times\Gamma_2)\ar[d]\ar[r]^-{(f_1\times f_2)_*}&\mathrm{H}_{i+j}(\Lambda\times\Lambda)\ar[d]\\
    &\mathrm{H}_{i+j}(U,\sim\Gamma_1\times_M \Gamma_2)\ar[d]^\cong\ar[r]^-{(f_1\times f_2)_*}&\mathrm{H}_{i+j}(U_{CS},\sim \Lambda\times_M \Lambda)\ar[d]^\cong\\
    \mathrm{H}_{i+j-\lambda_1-\lambda_2-n}(\Sigma\times_M \Sigma)\ar[d]_-{(\Delta')_!}\ar[r]^-{(p_1\times p_2)_!}\ar@{}[rd]|\dagger&\mathrm{H}_{i+j-n}(\Gamma_1\times_M \Gamma_2)\ar[d]^{\Phi_*}_\cong\ar[r]_-{(f_1\times f_2)_*}&\mathrm{H}_{i+j-n}(\Lambda\times_M \Lambda)\ar[d]^-{\mathrm{concat}_*}\\
    \mathrm{H}_{i+j-\lambda_1-\lambda_2-n-\mu}(\Sigma)\ar[r]_-{(p_3)_!}&\mathrm{H}_{i+j-n}(\Gamma_3)\ar[r]_-{(f_3)_*}&\mathrm{H}_{i+j-n}(\Lambda)
    }
\end{equation*}

Apart from the square at the bottom left corner, the above diagram is commutative, since the composition in the centre of the diagram, $\mathrm{H}_{*}(\Gamma_1\times\Gamma_2)\to\mathrm{H}_{*-n}(\Gamma_1\times_M\Gamma_2)$, is the Gysin map of the inclusion.
The square marked with $\dagger$ commutes only on the image of the injection $$\phi_!\colon \mathrm{H}_{*}(\Sigma)\to \mathrm{H}_{*+\mu}(\Sigma\times_M \Sigma)\,,$$ 
i.e. if $y\in \mathrm{H}_{*+\mu}(\Sigma\times_M \Sigma)$ is an element of the image of $\phi_!$, then 
$$\Phi_*\circ(p_1\times p_2)_!(y)=(p_3)_!\circ(\Delta')_!(y) \, .$$ This follows from the second commutative diagram of Lemma \ref{completing_manifolds_diffeo} and from $\Psi_!=\Phi_*$, where $\Psi$ denotes the inverse of $\Phi$ given in the proof of Lemma \ref{completing_manifolds_diffeo}.
Since $p_i$ has the section $s_i$ it follows that $(p_i)_!\colon \mathrm{H}_*(\Sigma)\to \mathrm{H}_{*+\lambda_i}(\Gamma_i)$ is injective. Using the sections $s_1\times s_2$ and $(s_1\times s_2)_{|_{\Sigma\times_M \Sigma}}$ the same is true for the maps 
$$(p_1\times p_2)_!\colon \mathrm{H}_*(\Sigma\times \Sigma)\to \mathrm{H}_*(\Gamma_1\times\Gamma_2) $$ and $$ (p_1\times p_2)_!\colon \mathrm{H}_*(\Sigma\times_M \Sigma)\to \mathrm{H}_*(\Gamma_1\times_M\Gamma_2) \,   . $$ It follows that all the horizontal compositions in the above diagram are injective, see Remark \ref{splitting_sublevels}.
We have that
\begin{equation*}
    \mathrm{H}_*(\Sigma\times_M\Sigma)\cong\mathrm{ker}\left((\Delta')_!\right)\oplus\mathrm{im}(\phi_!)\cong \mathrm{ker}\left((\Delta')_!\right)\oplus \mathrm{H}_{*-\mu}(\Sigma)\,.
\end{equation*}
Let $x+y=J_!(a\times b)$ denote the corresponding splitting of the image of $a\otimes b$ inside $\mathrm{H}_*(\Sigma\times_M\Sigma)$, i.e. $x\in\mathrm{ker}\left((\Delta')_!\right)$ and $y\in\mathrm{im}(\phi_!)$.
Then we have
\begin{align*}
(f_3)_*\circ(p_3)_!\circ(\Delta')_!(y)&=
   (f_3)_*\circ(p_3)_!\circ(\Delta')_!(x+y)\\&
   = (f_3)_*\circ(p_3)_!\circ(\Delta')_!\left(J_!(a\times b)\right)\\
    &=  (f_3)_*\circ(p_3)_!(a\bullet b)\neq0\,,
\end{align*}
which is nonzero since $(f_3)_*\circ(p_3)_!$ is injective and $a\bullet b\neq 0$ by assumption. On the other hand, by the commutativity of the rest of the diagram we have
\begin{align*}
    (f_1)_*\circ(p_1)_!(a)\wedge&(f_2)_*\circ(p_2)_!(b)=\mathrm{concat}_*\circ(f_1\times f_2)_*\circ(p_1\times p_2)_!(x+y)\\
    &=\mathrm{concat}_*\circ(f_1\times f_2)_*\circ(p_1\times p_2)_!(x)+\mathrm{concat}_*\circ(f_1\times f_2)_*\circ(p_1\times p_2)_!(y)\\
    &=(f_3)_*\circ\Phi_*\circ(p_1\times p_2)_!(x) + (f_3)_*\circ(p_3)_!(a\bullet b)\,.
\end{align*}
We now show that this last expression cannot be zero by showing that $(f_3)_*\circ\Phi_*\circ(p_1\times p_2)_!(x)$ cannot be a multiple of $ (f_3)_*\circ(p_3)_!(a\bullet b)$. One checks that the diagram
\begin{equation*}
    \xymatrix{
    \Gamma_1\times_M\Gamma_2\ar[d]_-{p_1\times p_2}&\Gamma_3\ar[l]_-{\Psi}\\
    \Sigma\times_M\Sigma&\Sigma\ar[u]_-{s_3}\ar[l]_-{\Delta'}
    }
\end{equation*}
commutes.
It follows that $$(s_3)_!\circ\Phi_*\circ(p_1\times p_2)_!=(\Delta')_! \, .$$
Hence, $\Phi_*\circ(p_1\times p_2)_!(x)$ is an element of the kernel of $(s_3)_!$, since $x\in \mathrm{ker}((\Delta')_!)$. 
This in turn means that $$  \Phi_*\circ(p_1\times p_2)_!(x) \in   \mathrm{H}_*(\Gamma_3\setminus \Sigma) \subseteq \mathrm{H}_*(\Gamma_3) \, .$$
The element $(p_3)_!(a\bullet b)$ lies in $\mathrm{H}_*(\Gamma_3,\Gamma_3\setminus\Sigma)\cong \mathrm{H}_{*-\mu}(\Sigma)\cong\mathrm{im}\left((p_3)_!\right)$, i.e. in the other direct summand of $\mathrm{H}_*(\Gamma_3)$, see Section \ref{sec_completing}. Let $L$ denote the energy of the closed geodesic $\gamma_3$.
 By Proposition \ref{prop_completing_mfld} we have the direct sum decomposition
\begin{equation*} 
        \mathrm{H}_{*}(\Lambda^{\leq L}) \cong \mathrm{H}_{*}(\Lambda^{<L}) \oplus \mathrm{H}_{*}(\Lambda^{\leq L},X^{<L})\,.
\end{equation*}
It follows from Remark \ref{splitting_sublevels} that under this decomposition $(f_3)_*\circ\Phi_*\circ(p_1\times p_2)_!(x)$ is an element of $\mathrm{H}_*(\Lambda^{<L})$ while $(f_3)_*\circ(p_3)_!(a\bullet b)$ lies in $\mathrm{H}_*(\Lambda^{\leq L},\Lambda^{<L})$. Consequently, the two terms cannot cancel.
\end{proof}

A version of the above theorem for a level Chas-Sullivan product in the case of rank one symmetric spaces is Theorem 13.5 in \cite{goresky:2009}.
The following corollary of the theorem is an extension of Proposition \ref{cs_prod_orientation_classes}:

\begin{cor}
For a homology class $a\in \mathrm{H}_*(\Sigma)$ of $\Sigma$ we have
\begin{equation*}
    (f_1)_*\circ(p_1)_!(a)\wedge(f_2)_*\circ(p_2)_!([\Sigma])=(f_3)_*\circ(p_3)_!(a)\,.
\end{equation*}
\end{cor}

\begin{proof}
We use the notation of the proof of the above theorem.
If $J_!(a\times b)=y\in\mathrm{im}(\phi_!)$ then we have
\begin{equation*}
    (f_1)_*\circ(p_1)_!(a)\wedge(f_2)_*\circ(p_2)_!(b)=(f_3)_*\circ(p_3)_!(a\bullet b)\,.
\end{equation*}
This is the case when $b = [\Sigma]$: the intersection product $a\bullet[\Sigma]$ of $[\Sigma]$ with any other class $a$ is given by $(\Delta')_!\circ J_!(a\times[\Sigma])$ and equals $a$. Since the projection onto the first factor $pr_1\colon\Sigma\times\Sigma\to\Sigma$ satisfies $\phi=pr_1\circ J$ and $(pr_1)_!(a)=a\times[\Sigma]$ we have
\begin{equation*}
    \phi_!(a)=J_!\circ (pr_1)_!(a)=J_!(a\times[\Sigma])\,.
\end{equation*}
As the square marked with $\dagger$ in the above diagram commutes on the image on $\phi_!$, the statement of the corollary follows.
\end{proof}

The case $a=[\Sigma]$ gives an alternative proof for Proposition \ref{cs_prod_orientation_classes} and Corollary \ref{cor_non_nilpotent}. Namely, as we have $[\Sigma]\bullet[\Sigma]=[\Sigma]$ it follows from Theorem \ref{theorem_cs_intersection} that the class
\begin{equation*}
    \Theta_\sigma= (f_\sigma)_*\circ (p_\sigma)_!\left([\Sigma]\right)\,,
\end{equation*}
which is associated to a prime closed geodesics $\sigma$, satisfies $\Theta_\sigma^{\wedge n}\neq0$ for all $n\in\mathbb{N}$. Here the superscript $\wedge n$ denotes the $n$-fold Chas-Sullivan product.

Combining this with Proposition \ref{infty_many_prime_geodesics} from Section \ref{sec_symm_geodesics}, which states that for symmetric spaces of higher rank there are infinitely many geometrically distinct prime critical submanifolds, we can summarize the findings of this section in the following theorem.

\begin{theorem}\label{nonnilpotent_CS_classes}
Let $M$ be a compact symmetric space and take homology with $\mathbb{Z}_2$-coefficients.
Every critical manifold of prime closed geodesics determines a non-nilpotent element in the Chas-Sullivan algebra of $\Lambda M$. If the rank of $M$ is two or larger, there are infinitely many such elements.
Such elements also exist if $M$ is only homotopy equivalent to a compact symmetric space, as the Chas-Sullivan product is a homotopy invariant.

Moreover, if $\sigma$ denotes a prime closed geodesic of $M$, then multiplication with $\Theta_\sigma$ gives an isomorphism from the loop space homology generated by $\sigma^m$ to the homology generated by $\sigma^{ m+1}$ for every $m\in\NN$.
More precisely, let $m\in\mathbb{N}$ and let $$U_{\sigma^m}  = \mathrm{im}\left((f_{\sigma^m})_*\circ (p_{\sigma^m})_! \colon \mathrm{H}_*(\Sigma) \hookrightarrow \mathrm{H}_{*+ \mathrm{ind}(\sigma^m)}(\Lambda M)  \right)$$ be the subspace of $\mathrm{H}_*(\Lambda M)$ induced by the critical manifold $\Sigma_{\sigma^m}$ that contains $\sigma^m$.
Then the map
\begin{equation*}
    \wedge \, \Theta_\sigma \colon U_{\sigma^m} \to U_{\sigma^{m+1}} \,,\,h\mapsto h\wedge\Theta_\sigma\,,
\end{equation*}
is an isomorphism for every $m\in\NN$.
\end{theorem}
\begin{remark}
    The above theorem is only stated for $\mathbb{Z}_2$-coefficients.
    However, we remark that if $R$ is a commutative ring with unit and all critical manifolds and completing manifolds are orientable, then the constructions of this section can be done for homology with coefficients in $R$.
    It remains to choose the orientations of the critical manifolds and the completing manifolds and figure out the signs.
    
    In particular if $M$ is a compact symmetric space and $\Sigma\subseteq \Lambda M$ is a critical manifold such that both $\Sigma$ and the corresponding completing manifold $\Gamma$ is orientable, then the homology class $f_*[\Gamma]\in \mathrm{H}_{\bullet}(\Lambda M)$ is a non-nilpotent class in the Chas-Sullivan ring of $M$.    
    Ziller discusses the orientability of the completing manifolds in \cite[p. 19-20] {ziller:1977}.
    In particular, the completing manifolds of symmetric spaces of higher rank are in general not orientable, even if the underlying manifold is simply connected.
\end{remark}

\section{Triviality of the based coproduct for symmetric spaces of higher rank} \label{sec_coproduct}

Using Bott's $K$-cycles in the based loop space of a compact symmetric space we will show in this section that the based coproduct is trivial for symmetric spaces of higher rank.
We will then study the implications of this result for the string topology coproduct on the free loop space.
We start by defining Bott's $K$-cycles in symmetric spaces following \cite{bott:1958a}.
See also \cite{araki:1962} for more properties of these manifolds.

Let $M = G/K$ be a compact simply connected symmetric space.
We adopt all constructions and notations from Section \ref{sec_symm_geodesics}.
Consider a pair $(\alpha,n)$ where $\alpha\in \Delta$ is a root and $n\in \ZZ$.
As seen in Section \ref{sec_symm_geodesics} this defines an affine hyperplane
$$   \{H\in \mathfrak{a}\,|\, \alpha(H) = n \}\subseteq \mathfrak{a} \,  ,   $$
the singular plane $(\alpha,n)$.
Denote the stabilizer of the maximal torus $T\subseteq M$ under the action of $K$ on $M$ by $K_0$.
If $p=(\alpha,n)$ is a singular plane, then let $\overline{p} = \mathrm{Exp}(p)\subseteq T$ and set
$$   K_p \coloneqq \{k\in K\,|\, k.q = q \,\, \text{for all} \,\, q\in \overline{p} \} \, .    $$
It is clear that both $K_p$ and $K_0$ are closed subgroups of $K$ and that $K_0\subseteq K_p$.
Furthermore one can show that
$$  \mathrm{dim}(K_p) > \mathrm{dim}(K_0)\, , $$
for each singular plane $p$, see \cite[p. 91]{araki:1962}.

Let $P = (p_1,\ldots,p_m)$ be an ordered family of singular planes.
We set
$$  W(P)  \coloneqq K_{p_1}\times K_{p_2}   \times \ldots \times K_{p_m}\,  . $$
There is a right-action of the $m$-fold product of $K_0$ on $W(P)$, that is given by 
\begin{eqnarray*}    \chi : W(P) \times K_0^m &\to& W(P)      \\
                ( (k_1,\ldots,k_m), (h_1,\ldots,h_m)) &\mapsto& (k_1 h_1,h_1^{-1} k_2 h_2,\ldots, h_{m-1}^{-1} k_m h_m)\,   .
\end{eqnarray*} 
This is clearly a proper and free action and consequently the quotient
$$  \Gamma_P = W(P)/ (K_0^m)    $$
is a compact smooth manifold. 
Note the similarity of the construction of these manifolds with the construction of the completing manifolds in Section \ref{sec_chas_sullivan}.

We now define an embedding of $\Gamma_P$ into the based loop space $\Omega_o M$.
Let $c^P = (c_{0}^P,c_{1}^P,\ldots,c_{m}^P)$ be an ordered family of polygons in $\mathfrak{a}$ such that the following properties are satisfied:
\begin{itemize}
    \item The polygon $c_{0}^P$ starts at the origin $0\in\mathfrak{a}$.
    \item For $i=1,\ldots,m$ the endpoint of the polygon $c_{i-1}^P$ is the start point of the polygon $c_{i}^P$ and lies on the singular plane $p_i$.
    \item The polygon $c_{m}^P$ ends at $0\in\mathfrak{a}$.
\end{itemize}
For $i=0,1,\ldots,m$, we parametrize the polygon $c_{i}^P$ on the interval $[\frac{i}{m+1},\frac{i+1}{m+1}]$.
Define a map $\widetilde{f}_P : W(P) \to \Omega_o M$ by setting
$$   \widetilde{f}_P(k_1,k_2,\ldots,k_m) \coloneqq \mathrm{concat}\Big[  \mathrm{Exp}(c_{0}^P), k_1\mathrm{Exp}(c_{1}^P), \ldots,k_1\ldots k_m\mathrm{Exp}(c_{m}^P)    \Big]\, ,     $$
where $\mathrm{concat}$ is the obvious concatenation map of multiple paths.
By the particular choice of the polygons, the map $\widetilde{f}_P$ is well-defined and continuous. 
The map $\widetilde{f}_P$ is invariant under the right-action $\chi$ of $K_0^m$.
Consequently it induces a map
$$  f_P : \Gamma_P \to \Omega_o M     $$
and one can check that this is indeed an embedding.
Note that the homotopy type of $f_P$ does not depend on the specific choice of polygons $c^P$, see the discussion before Definition 5.4 in \cite{stegemeyer:2021}.

Let $R$ be a commutative ring with unit.
If the manifold $\Gamma_P$ is orientable with respect to $R$ and if $[\Gamma_P]\in \mathrm{H}_{*}(\Gamma_P;R)$ denotes a corresponding fundamental class then we obtain a homology class
$$   P_* \coloneqq (f_P)_*[\Gamma_P] \in \mathrm{H}_{*}(\Omega_o M; R)\,  .  $$
If each $\Gamma_P$ is orientable with respect to $R$, Bott and Samelson show that the set
$$  \{ P_*\in \mathrm{H}_{*}(\Omega_o M)\,|\, P \,\,\text{is an ordered family of singular planes}\}     $$
generates the homology $\mathrm{H}_{*}(\Omega_o M)$, see \cite{bott:1958a}.
In particular, for every compact symmetric space we get a set that generates the homology of the based loop space with $\ZZ_2$-coefficients.
We will comment on the orientability assumption later on.
Note that the empty set is also seen as an ordered family of singular planes and the corresponding homology class corresponds to the basepoint $[o]\in \mathrm{H}_0(\Omega M)$.
\begin{remark}
In the present section we do not refer to the manifolds $\Gamma_P$ as completing manifolds, since we use the mentioned results by Bott and Samelson that the classes $P_*$ generate the homology of the based loop space.

However, we briefly want to sketch that the manifolds $\Gamma_P$ arise as completing manifolds in the following setting.
If $M$ is a connected manifold then it is well-known that the based loop space $\Omega M$ is homotopy equivalent to the space $\Omega_{pq}M$ of paths connecting two arbitrary points $p,q\in M$.
In the case of symmetric spaces, we choose $p=o$ to be the basepoint and we can choose the other point $q$ to lie in a given maximal torus $T$.
If $q$ is not conjugate to $o$ along any geodesic, then one can show
that all geodesics joining $o$ and $q$ are of the form $\gamma_H$, where $H\in\mathfrak{a}$, as in Section \ref{sec_symm_geodesics}.
Recall that $\gamma_H = \Exp \circ \, \sigma_H$ with $\sigma_H$ being a straight line in $\mathfrak{a}$ connecting $0$ and $H$.
Moreover, we have seen in Section \ref{sec_symm_geodesics} that the conjugate points along $\gamma_H$ appear at times when $\sigma_H$ intersects a singular plane.
Now, one can argue that if $q$ is in \textit{general position} then $\sigma_H$ never meets an intersection of several hyperplanes, see \cite{bott:1958a}.
Hence, if $P = (p_1,\ldots,p_m)$ is the ordered sequence of singular planes that $\sigma_H$ intersects, then up to the last part of the polygon that goes back to the origin, we can choose the individual bits of $\sigma_H$ for the construction of the sequence of polygons $(c_0^P,\ldots,c_m^P)$.

Therefore, if we do Morse theory in the path space $\Omega_{oq} M$ then for appropriate families of singular planes $P$ the manifolds $\Gamma_P$ actually become completing manifolds for the critical points of the energy functional, which are just the geodesics connecting $o$ and $q$.
For details, see \cite{bott:1958a}.
\end{remark}

We will now show that the basepoint intersection multiplicity of a class $P_*$, where $P$ is an ordered family of singular planes, is $1$.
As in Section \ref{sec_symm_geodesics} let $\mathcal{F}\subseteq \mathfrak{a}$ be the lattice
$$ \mathcal{F}  = \{ H\in \mathfrak{a}\,|\, \mathrm{Exp}(H) = o\} \, .    $$
We make the following definition.
\begin{definition}
Let $P = (p_1,\ldots,p_m)$ be an ordered family of singular planes.
An ordered family of polygons $(c_0^P,c_1^P,\ldots,c_m^P)$ is called \textit{lattice-nonintersecting} if the following holds: No polygon $c_i^P$ intersects the lattice $\mathcal{F}$ apart from $c_0^P$ at its start point and $c_m^P$ at its endpoint.
\end{definition}
If the rank of $M$ is greater or equal than $2$, i.e. $\mathrm{dim}(\mathfrak{a})\geq 2$, then it is clear that for any non-trivial ordered family of singular planes $P$ we can choose a corresponding ordered family of polygons which is lattice-nonintersecting. Recall that the class $P_*\in \mathrm{H}_{*}(\Omega M)$ does not depend on this choice.
\begin{prop} \label{prop_bott_intersection_m}
Let $M$ be a compact simply connected symmetric space of rank $r\geq 2$ and let $P$ be an ordered family of singular planes.
\begin{enumerate}
    \item The induced class $P_*\in \mathrm{H}_{*}(\Omega M,o;\ZZ_2)$ satisfies $\mathrm{int}(P_*) = 1$.
    \item If $\Gamma_P$ is orientable the same holds with arbitrary coefficients. 
\end{enumerate}
\end{prop}
\begin{proof}
Choose a representing cycle $X$ of $f_*[\Gamma_P]$ and let $\gamma\in \mathrm{im}(X)$.
Since the rank of $M$ is greater or equal than $2$, we can define the class $P_*$ via a lattice-nonintersecting family of polygons $(c_0^P,\ldots,c_m^P)$.
Assume that there is $s\in(0,1)$ with $\gamma(s) = o$.
Then there are $(k_1,\ldots,k_m)\in W(P)$ with
$$  k_1\ldots k_i \mathrm{Exp}(c_i^P(s)) = o     $$
for some $i\in\{1,\ldots,m\}$.
But then $\mathrm{Exp}(c^P_i(s)) = o$, so $c_i^P(s)\in \mathcal{F}$, which contradicts our choice of a lattice-nonintersecting family of polygons.
\end{proof}
\begin{cor} \label{cor_triviality_based}
Let $M$ be a compact simply connected symmetric space of rank $r\geq 2$.
\begin{enumerate}
    \item The based string topology coproduct $\cpr_{\Omega}$ on $\mathrm{H}_{*}(\Omega M,o;\ZZ_2)$ is trivial.
\item If the manifold $\Gamma_P$ is orientable for each ordered family of singular planes $P$, then the based string topology coproduct is trivial for homology with arbitrary coefficients.
\end{enumerate}
\end{cor}
\begin{proof}
This follows from Propositions \ref{prop_intersection_trivial} and \ref{prop_bott_intersection_m} and the fact that the classes of the form $P_*$ with $P$ being a non-trivial ordered family of singular planes generate the homology $\mathrm{H}_*(\Omega_o M,o)$.
\end{proof}
\begin{remark}
Araki determines the irreducible symmetric spaces for which the orientability assumption of Corollary \ref{cor_triviality_based}.(2) holds, see \cite[Corollary 5.12]{araki:1962}.
These include all compact Lie groups, the complex Grassmannians, the quaternionic Grassmannians, spheres, the real Grassmannians of the form $SO(2n+2)/SO(2)\times SO(n)$, the manifolds of the form $SU(2n)/Sp(n)$ and of the form $SO(2n)/U(n)$ as well as some quotients of exceptional Lie groups.
See \cite[Corollary 5.12]{araki:1962} for the complete list.
\end{remark}

It is now natural to ask whether the string topology coproduct on the free loop space of compact symmetric spaces of higher rank is trivial as well.
Although we cannot answer this question, we will show in the following that the triviality of the based coproduct implies the triviality of the coproduct on a certain subspace of the homology of the free loop space.
As seen in Section \ref{sec_def} the diagram
$$
\begin{tikzcd}
    \mathrm{H}_i(\Omega M,o) \arrow[]{r}{i_*} \arrow[]{d}{\vee_{\Omega}} & \mathrm{H}_i(\Lambda ,M) \arrow[]{d}{\vee}
    \\
    \mathrm{H}_{i+1-n}(\Omega \times \Omega,\Omega\times o\cup o \times \Omega) \arrow[]{r}{(i\times i)_*} & \mathrm{H}_{i+1-n}(\Lambda \times \Lambda,\Lambda\times M\cup M\times \Lambda) .
\end{tikzcd}
$$
commutes.
Hence, it is clear that $\vee$ is trivial on the image of $i_*\colon \mathrm{H}_{*}(\Omega M,o)\to \mathrm{H}_{*}(\Lambda ,M)$ if $\vee_{\Omega}$ is already trivial. 
Moreover, note that any action $\psi\colon H\times M\to M$ of a compact Lie group $H$ induces an action $\Psi\colon H\times \Lambda M\to \Lambda M$ on the free loop space, which preserves the constant loops.
In particular, we obtain an induced pairing in homology
$$      \Theta\colon \mathrm{H}_i(H)\otimes \mathrm{H}_j(\Lambda M,M)  \to \mathrm{H}_{i+j}(\Lambda M,M),\quad \Theta(X\otimes Y) = \Psi_* ( X\times Y),  $$
where $ X\in\mathrm{H}_i(H), Y\in \mathrm{H}_j(\Lambda M,M)$ .  
\begin{prop} \label{prop_triviality_in_free_loop_space}
    Let $M$ be a compact simply connected symmetric space of rank greater than or equal to $2$.
    Consider homology with coefficients in an arbitrary commutative unital ring $R$.
    Assume that $H$ is a compact Lie group acting on $M$ and consider the pairing $\Theta$ as above.
    Let $\mathcal{K}$ be the image of $\Theta$ when restricted to the subspace
    $$     \mathrm{H}_{*}(H) \otimes \big(i_* \,(\mathrm{H}_{*}(\Omega M,o))\big) \subseteq \mathrm{H}_{*}(H)\otimes \mathrm{H}_{*}(\Lambda M,M) .       $$
    Then the string topology coproduct is trivial on $\mathcal{K}$.
    In particular $i_*(\mathrm{H}_*(\Omega M,o))$ is contained in $\mathcal{K}$.
\end{prop}
\begin{proof}
    Let $Z = \Theta(X\otimes Y)$ be in $\mathcal{K}$, i.e. $X\in\mathrm{H}_{*}(H)$, $Y = i_* W$ for $W\in\mathrm{H}_*(\Omega M,o)$.
    Assume that $W$ and $i_* W$ are non-trivial.
    As we have seen previously, we can choose a cycle representing $W$ with no non-trivial basepoint intersection in $\Omega M$.
    But then it is clear directly that $i_* W$ has a representing cycle which does not have any non-trivial basepoint intersections in $\Lambda M$ either.
    Moreover, since a group acts via diffeomorphisms, the group action does not increase basepoint intersection multiplicity, hence we see that every homology class in $\mathcal{K}$ has basepoint intersection multiplicity equal to $1$.
    By \cite[Theorem 3.10]{hingston:2017} this implies $\vee Z = 0$.
\end{proof}
\begin{remark}
    Note that in general it is not clear \textit{how much} of $\mathrm{H}_*(\Lambda M,M)$ we obtain by considering group actions as above.
    Hepworth computes the pairing $\Theta$ in the case of $M = \mathbb{C}P^n$ and $H = U(n+1)$.
    In this case $\mathcal{K}$ gives a considerably larger subspace of $\mathrm{H}_{*}(\Lambda M,M)$ than just $i_* (\mathrm{H}_*(\Omega M,o))$.
    Moreover, note that in the case of $M$ being a compact Lie group, the subspace $\mathcal{K}$ is of course all of $\mathrm{H}_*(\Lambda M,M)$ as we have $\Lambda M\cong \Omega M\times M$ in this case, see e.g. \cite{stegemeyer:2021}.
    Therefore we believe that for an arbitrary compact symmetric space of higher rank the subspace $\mathcal{K}$ will contain $i_*(\mathrm{H}_*(\Omega M,o))$ as a proper subspace.
\end{remark}

To conclude this section we show that for a compact symmetric space $M$ we can associate a homology class $[W_{\Sigma}]\in\mathrm{H}_*(\Lambda M,M;\ZZ_2)$, which lies in $i_*(\mathrm{H}_*(\Omega M,o;\ZZ_2))$, to any critical manifold $\Sigma\subseteq \Lambda M$ .
This will be done by using a version of Ziller's completing manifolds on the based loop space.
For the rest of this section we consider homology and cohomology with $\ZZ_2$-coefficients.

In the following let $M =G/K$ be a compact simply connected symmetric space of arbitrary rank.
We saw that the critical manifolds of the energy functional in the free loop space are of the form $\Sigma = G/K_c$, where $c$ is a closed geodesic.
Correspondingly, the critical manifolds in the based loop space are of the form $\Sigma_{\Omega} = K/K_c$. We can build a completing manifold $\Gamma_{\Omega}$ in the based loop space just as we did in the free loop space.
More explicitly, we set
$$ \Gamma_{\Omega}    =  (K\times K_1\times \ldots K_l)/ (K_c)^{l+1}    $$
where the groups $K_i$ are the stabilizers of the conjugate points along $c$ exactly as we had it in the case of the free loop space, see Section \ref{sec_chas_sullivan}.
Moreover, note that the action of the group $(K_c)^{l+1}$ on $K\times K_1\ldots\times K_l$ is obtained by restricting the action $\chi$, which was defined in equation \eqref{eq_definition_chi}.
If $\Gamma$ is the completing manifold associated to $\Sigma$ then we can embed $\Gamma_{\Omega}$ into $\Gamma$.
Moreover, the embedding $f\colon \Gamma\to \Lambda M$ restricts to an embedding of $\Gamma_{\Omega}$ into the based loop space.
The submersion $p\colon \Gamma\to \Sigma$ restricts to a submersion $p_{\Omega}\colon \Gamma_{\Omega}\to \Sigma_{\Omega}$.
These compatabilities are summarized in the following two commutative diagrams.
$$   
\begin{tikzcd}
    \Gamma_{\Omega} \arrow[]{r}{f|_{\Gamma_{\Omega}}} \arrow[hook, swap]{d}{i_{\Gamma}}
    &
    \Omega M \arrow[hook]{d}{i} \\
    \Gamma \arrow[]{r}{f} & \Lambda 
\end{tikzcd}
\qquad \qquad
\begin{tikzcd}
    \Gamma_{\Omega} \arrow[hook]{r}{i_{\Gamma}}  \arrow[swap]{d}{p_{\Omega}} 
    & \Gamma \arrow[]{d}{p}
    \\
    \Sigma_{\Omega} \arrow[hook]{r}{i_{\Sigma}}
    & \Sigma
\end{tikzcd}
$$
Moreover, the section $s\colon \Sigma\to \Gamma$ induces a section $s_{\Omega}\colon \Sigma_{\Omega}\to\Gamma_{\Omega}$ by restriction.
All these properties can be checked by using the explicit expressions of the respective maps, which were introduced in Section \ref{sec_chas_sullivan}.

\begin{lemma} \label{lemma_compatibility_completing_manifolds}
The following diagram commutes up to sign.
$$\begin{tikzcd}
    \mathrm{H}_i(\Sigma_{\Omega}) \arrow[]{r}{(p_{\Omega})_!} \arrow[swap]{d}{(i_{\Omega})_*} 
    & \mathrm{H}_{i+\lambda}(\Gamma_{\Omega}) \arrow[]{d}{(i_{\Gamma})_*}
    \\
    \mathrm{H}_i(\Sigma) \arrow[]{r}{p_!} & \mathrm{H}_{i+\lambda}(\Gamma)
\end{tikzcd}$$
where $i_{\Gamma}$ and $i_{\Omega}$ are the obvious inclusions.
\end{lemma}
\begin{proof}
First, consider the diagram
\begin{equation} \label{diagram_thom_omega}
\begin{tikzcd}
    \mathrm{H}_i(\Sigma_{\Omega}) \arrow[<-]{r}{(s_{\Omega})_!} \arrow[swap]{d}{(i_{\Omega})_*} 
    & \mathrm{H}_{i+\lambda}(\Gamma_{\Omega}) \arrow[]{d}{(i_{\Gamma})_*}
    \\
    \mathrm{H}_i(\Sigma) \arrow[<-]{r}{s_!} & \mathrm{H}_{i+\lambda}(\Gamma)
\end{tikzcd}
\end{equation}
and recall that $s_!$ is equal to the composition
$$   \mathrm{H}_{i}(\Gamma) \to \mathrm{H}_i(\Gamma,\Gamma\setminus \Sigma) \xrightarrow[]{\text{excision}} \mathrm{H}_i(U,U\setminus \Sigma) \xrightarrow[]{\text{Thom}} \mathrm{H}_{i-\lambda}(\Sigma)   $$
up to sign. 
Here, $U$ is a tubular neighborhood of $\Sigma$ in $\Gamma$.
The codimension of $\Sigma_{\Omega}$ in $\Gamma_{\Omega}$ is the same as the one of $\Sigma$ in $\Gamma$.
The Thom class of $(U,U\setminus \Sigma)$ pulls back to the Thom class of $(U_{\Omega},U_{\Omega}\setminus \Sigma_{\Omega})$ where $U_{\Omega}$ is an appropriatley chosen tubular neighborhood of $\Sigma_{\Omega}$ in $\Gamma_{\Omega}$.
Consequently, the above diagram commutes.
Now, let $X\in \mathrm{H}_*(\Sigma_{\Omega})$, then we have
$$    X =  (p_{\Omega}\circ s_{\Omega})_! X = (s_{\Omega})_! (p_{\Omega})_! X , $$
since $p_{\Omega}\circ s_{\Omega} = \id_{\Sigma_{\Omega}}$.
Thus, we get
$$     p_! (i_{\Omega})_* X = p_! (i_{\Omega})_* (s_{\Omega})_! (p_{\Omega})_! X = \pm p_!  s_! (i_{\Gamma})_* (p_{\Omega})_! X = \pm (i_{\Gamma})_* (p_{\Omega})_! X    $$
where we used the commutativity of the diagram \eqref{diagram_thom_omega}.
This completes the proof.
\end{proof}
\begin{prop} \label{prop_wc_classes}
    Let $M$ be a compact simply connected symmetric space of rank greater than or equal to $2$.
    To every critical submanifold $\Sigma$ we can associate a class $[W_{\Sigma}]\in \mathrm{H}_{\mathrm{ind}(\Sigma)}(\Lambda M,M; \ZZ_2)$ that has trivial coproduct.
\end{prop}
\begin{proof}
Let $\Sigma$ be a critical manifold and denote by $\Gamma$ its associated completing manifold, which comes as usual with the maps $p$ and $f$.
    We define $$[W_{\Sigma}] =  f_*\circ p_! [x_0]  $$
    where $[x_0]\in\mathrm{H}_0(\Sigma)$ is a generator in degree $0$.
    Since $\Gamma$ is a completing manifold, we have $[W_{\Sigma}]\neq 0$.
    By Lemma \ref{lemma_compatibility_completing_manifolds} we see that the following diagram commutes.
    $$
\begin{tikzcd}
     \mathrm{H}_i(\Sigma_{\Omega}) \arrow[]{r}{(p_{\Omega})_!} \arrow[swap]{d}{(i_{\Omega})_*} 
    & \mathrm{H}_{i+\lambda}(\Gamma_{\Omega}) \arrow[]{d}{(i_{\Gamma})_*} \arrow[]{r}{(f_{\Gamma_{\Omega}})_*}
& \mathrm{H}_{i+\lambda}(\Omega)     \arrow[]{d}{i_*}
    \\
    \mathrm{H}_i(\Sigma) \arrow[]{r}{p_!} & \mathrm{H}_{i+\lambda}(\Gamma) \arrow[]{r}{f_*}
    & \mathrm{H}_{i+\lambda}(\Lambda)
\end{tikzcd}
$$
Clearly the generator $[y_0]\in\mathrm{H}_0(\Sigma_{\Omega})$ in degree $0$ is mapped to $[x_0]$ under $(i_{\Omega})_*$.
If we follow the class $[y_0]$ around the diagram we see that
$$      [W_{\Sigma}] = i_* \, (f_{\Omega})_* \, (p_{\Omega})_! [y_0],     $$
so in particular $[W_{\Sigma}]\in i_*(\mathrm{H}_*(\Omega))$.
By Proposition \ref{prop_triviality_in_free_loop_space} we then have $\vee[W_{\Sigma}] = 0$.
\end{proof}
\begin{remark}\label{remark_index-growth} Let $M$ be a compact symmetric space of dimension $\mathrm{dim}(M) = n$.
\begin{enumerate}
    \item 
    We want to stress that in the rank $1$ case the classes of the form $[W_{\Sigma}]$ have non-trivial coproduct and that in addition their dual classes in cohomology are non-nilpotent elements for the Goresky-Hingston product.
    Recall that in Section \ref{sec_chas_sullivan} we observed that the classes $f_*[\Gamma]$ are non-nilpotent elements in the Chas-Sullivan product for any rank.
    In fact our study of these classes was motivated by the well-known behavior of these classes in the rank $1$ case.
    Thus, it would be natural to believe that the string topology coproduct of the classes of the form $[W_{\Sigma}]$ behaves analogously as well when comparing the rank $1$ and the higher rank case. 
    However, we have just seen that this is not the case.
    \item We further want to note that from a geometric point of view the difference between the Chas-Sullivan product and the Goresky-Hingston product that occurs for higher rank symmetric spaces seems to be related to the following observation.
    If $M$ is a compact rank one symmetric space, then the cohomology class $\omega\in\mathrm{H}^{*}(\Lambda M,M)$ which is dual to the homology class $[W_{{\Sigma}_{\gamma}}]$ is non-nilpotent for the Goresky-Hingston product.
    Here, $\Sigma_{\gamma}$ is the critical manifold which contains the prime closed geodesic $\gamma$.
    Since $M$ is of rank $1$, $\Sigma_{\gamma}$ thus contains all prime closed geodesics.
    Moreover, the iterates of $\omega$ under the Goresky-Hingston product correspond to the iteration of closed geodesics, i.e. the cohomology class $\omega^{\ostar k}\in\mathrm{H}^*(\Lambda M,M)$ is dual to the homology class $[W_{\Sigma_{\gamma^k}}]$.
    In particular, we have
    $$   \mathrm{deg}(\omega^{\ostar k}) = \mathrm{ind}(\gamma^k) ,     $$
    i.e. the degree growth of the powers of $\omega$ equals the index growth of the closed geodesics.
    As we have seen in Section \ref{sec_symm_geodesics} for an arbitrary compact symmetric space $M = G/K$ the index of the iterates of a prime closed geodesic $\gamma $ satisfies
    $$    \mathrm{ind}(\gamma^k) = k\, \mathrm{ind}(\gamma ) + (k-1) \mu  .  $$
    Here, $\mu = \mathrm{dim}(K\cdot \gamma) =  \mathrm{dim}(\Sigma_{\gamma} \cap \Omega M)$ is the dimension of the critical manifold containing $\gamma$ in the based loop space.
    Hence, we have
    $$    \mathrm{ind}(\gamma^{k+\ell}) = \mathrm{ind}(\gamma^k) + \mathrm{ind}(\gamma^{\ell}) + \mu .  $$
    As discussed at the end of Section \ref{sec_symm_geodesics}, it holds that $\mu = n-1$ if and only if the rank of $M$ is $1$ and $\mu < n-1$ otherwise.
    But it happens that $n-1$ is exactly the degree shift of the Goresky-Hingston product!
    Heuristically speaking, it is thus not surprising that a connection between iterates of closed geodesics and powers of cohomology classes with respect to the Goresky-Hingston product only holds in the rank $1$ case.

\noindent
    In contrast, note that for the Chas-Sullivan product the non-nilpotent classes $[\Gamma_{\Sigma_{\gamma}}]$ associated to a critical manifold $\Sigma_{\gamma}$ have degree $\mathrm{deg}([\Gamma_{\Sigma_{\gamma}}] )= (\mathrm{ind} + \mathrm{null})(\gamma) $.
    By the fact that $\mathrm{null}(\gamma) = n+\mu$, see Section \ref{sec_symm_geodesics}, and the expression for the index, we see that 
    \begin{equation}\label{eq_12}
          (\mathrm{ind} + \mathrm{null}) (\gamma^{k + \ell}) =  (\mathrm{ind} + \mathrm{null})(\gamma^k)  + (\mathrm{ind}+\mathrm{null})(\gamma^{\ell}) - n . 
    \end{equation} 
    This holds for all compact symmetric spaces, independently of the rank.
    Hence, there is no difference here between rank $1$ symmetric spaces and symmetric spaces of higher rank and the degree shift of the Chas-Sullivan product matches with the shift in the addition formula for index plus nullity \eqref{eq_12}.
\end{enumerate}
\end{remark}


\section{The string topology coproduct on products of symmetric spaces}\label{sec_spheres}

In this final section we study the string topology coproduct on the free loop space of the product of two compact symmetric spaces
by studying how Ziller's cycles behave on product spaces.

Let $M_1$ and $M_2$ be two compact symmetric spaces and consider $M = M_1\times M_2$ equipped with the product metric.
This is again a symmetric space.
Indeed, if $s_p\colon M_1\to M_1$ and $s_q\colon M_2\to M_2$ are geodesic involutions at $p\in M_1$ and $q\in M_2$,respectively, then the map $s_p\times s_q\colon M_1\times M_2\to M_1\times M_2$ is a geodesic involution at $(p,q)\in M_1\times M_2$ with respect to the product metric.
We want to show that any intrinsic Ziller cycle in $M$ is a product of Ziller cycles of $M_1$ and $M_2$.
Note that we have $M_1 \cong G_1/K_1$ and $M_2 \cong G_2/K_2$ and consequently $M \cong (G_1\times G_2)/(K_1\times K_2)$.
Clearly a loop $c = (c_1,c_2)$ in $M$ is a closed geodesic in $M$ if and only if $c_1$ is a closed geodesic in $M_1$ and $c_2$ is a closed geodesic in $M_2$.
Consequently, the critical sets in $\Lambda M$ are of the form
$   \Sigma_1\times \Sigma_2,     $
where $\Sigma_1$ runs over all critical manifolds in $\Lambda M_1$ including the trivial loops and $\Sigma_2$ runs over all critical manifolds in $\Lambda M_2$ including the trivial loops.
The product
$    M_1\times M_2   $
of the critical manifolds of trivial loops are precisely the trivial loops in $\Lambda M$.
The stabilizer of a closed geodesic $c= (c_1,c_2)$ is
$     K_c = K_{c_1} \times K_{c_2} $
where $K_{c_1}$ is the stabilizer of $c_1$ and $K_{c_2}$ is the stabilizer of $c_2$.
Note that if $t_*\in (0,1)$ is such that $c_1(t_*)$ is a conjugate point along $c_1$ and $c_2(t_*)$ is not a conjugate along $c_2$, then the closed geodesic does have a conjugate point at time $t=t_*$.
Indeed, the stabilizer of this conjugate point is $K_{c(t_*)}\times K_{c_2}$ which has $K_{c_1}\times K_{c_2}$ as a subgroup of smaller dimension.
This shows that $c(t_*)$ is a conjugate point, see also Section \ref{sec_symm_geodesics}.
Of course, if both $c_1$ and $c_2$ have a conjugate point at a time $t=t_*$ then it is clear as well that $c$ has a conjugate point at this time.

Let $\Gamma_1$ and $\Gamma_2$ be the Ziller cycles associated to the critical manifolds $\Sigma_1$ and $\Sigma_2$ that contain the closed geodesics $c_1$ and $c_2$, respectively.
Of course in $\Lambda M$ there is an intrinsically defined Ziller cycle $\Gamma$ for the closed geodesic $c = (c_1,c_2)$.
We want to study how $\Gamma$ is related to $\Gamma_1$ and $\Gamma_2$.
Let 
$$    0< t_1 < t_2 < \ldots < t_l < 1     $$
be the times such that $o\in M$ has a conjugate point along $c$.
As we have seen above if $t\in\{t_1,\ldots,t_l\}$ then $c_1$ or $c_2$ have a conjugate point.
Let $K_1',\ldots,K_l'\subseteq K_1\times K_2$ be the stabilizer groups of $c(t_1),\ldots,c(t_l)$, respectively.
By the above it is clear that these groups always split
$$   K_i' =    K_{1,i}\times K_{2,i} \quad \text{with}\quad K_{c_1}\subseteq K_{1,i}\subseteq K_1 \,\,\,\text{and}\,\,\, K_{c_2}\subseteq K_{2,i}\subseteq K_2 .       $$
As usual we define
$$     W(c) = (G_1\times G_2)\times K_1'\times \ldots K_l'          $$
and we consider the action $\chi\colon W(c)\times (K_c)^{l+1}\to W(c)$ as in Section \ref{sec_chas_sullivan}.
Note that this action splits as follows.
Define 
$$   W'(c_1) = G_1 \times K_{1,1} \times K_{1,2}\times \ldots \times K_{l,1}       \quad \text{and}\quad      W'(c_2) = G_2 \times K_{2,1} \times K_{2,2}\times \ldots \times K_{l,2} .      $$
There is an obvious $(K_{c_i})^{l+1}$ action on $W'(c_i)$, $i\in\{1,2\}$, which we denote by $\chi'_i$.
Moreover, we have obvious isomorphisms of Lie groups
$$     \Phi\colon W(c) \to W'(c_1) \times W'(c_2)\quad \text{and} \quad  \Psi\colon (K_c)^{l+1} \to (K_{c_1})^{l+1}\times (K_{c_2})^{l+1}  .    $$
It is easy to check that $\Phi$ is equivariant with respect to $\Psi$, i.e. the following diagram commutes.
$$
\begin{tikzcd}
W(c)\times (K_c)^{l+1} \arrow[]{r}{\chi}  \arrow[swap]{d}{T\circ (\Phi\times \Psi)}  &  [3em] 
W(c) \arrow[]{d}{\Phi}
\\
W'(c_1)\times  (K_{c_1})^{l+1}\times W'(c_2) \times (K_{c_2})^{l+1} \arrow[]{r}{\chi'_1\times \chi'_2} & 
W'(c_1)\times W'(c_2) .
\end{tikzcd}
$$
Here $T$ is the obvious swapping map.
Consequently, there are induced diffeomorphisms
$$   \Gamma_c \xrightarrow[]{\cong}   \Gamma'_{c_1}\times \Gamma'_{c_2}     \quad \text{where} \quad    \Gamma'_{c_i} =        W'(c_i) / (K_{c_i})^{l+1}     , \,\,\,i\in\{1,2\}  . $$
We now want to compare $\Gamma'_{c_1}$ with the Ziller cycle $\Gamma_{c_1}$, which is intrinsically defined in $M_1$.
Recall that $\Gamma_{c_1}$ is defined as
$$    \Gamma_{c_1} = W(c_1)/(K_{c_1})^{r+1}    \quad \text{where} \quad     W(c_1) = \widetilde{K}_1\times \ldots \times \widetilde{K}_r    $$
and the $\widetilde{K}_i$ are the stabilizers of the conjugate points along $c_1$.
Note that among the groups $K_{1,i}$ we have $K_{1,i} = K_{c_1}$ whenever $t_i$ is such that $c_2(t_i)$ is conjugate but $c_1(t_i)$ is not.
However, every group in the ordered family $(K_{1,1},\ldots,K_{1,l})$, which is not $K_{c_1}$ appears as a member of the ordered family $(\widetilde{K}_1,\ldots,\widetilde{K}_r)$.

\begin{lemma} \label{lemma_diffeo_gamma_prime}
The manifolds $\Gamma'_{c_1}$ and $\Gamma_{c_1}$ are diffeomorphic.
\end{lemma}
\begin{proof}
In order to write down the diffeomorphisms, we need to introduce some notation.
Let $1\leq s_1 < \ldots < s_r \leq l$ be the indices such that $K_{1,s_i}$ is not equal to $K_{c_1}$.
In this case we have $K_{1,s_i}\cong \widetilde{K}_i$, $i\in\{1,\ldots,r\}$.
We define $\widetilde{\xi} \colon W_{c_1}\to W'_{c_1}$ and $\widetilde{\zeta} \colon W'_{c_1}\to W_{c_1}$ as  
$$      \widetilde{\xi} (\widetilde{k}_1,\ldots, \widetilde{k}_r) = (h_i)_{i\in\{1,\ldots,l\} } \quad \text{with}\quad  h_i = \begin{cases}  e,  & i \not\in \{s_1,\ldots,s_r\} \\ \widetilde{k}_j ,   & i = s_j .    \end{cases}       $$
and
$$    \widetilde{\zeta}( h_1,\ldots, h_{s_1-1},h_{s_1},\ldots, h_l) = (h_1\cdot h_2\cdot \ldots \cdot h_{s_1}, h_{s_1+1}\cdot \ldots \cdot h_{s_2},\ldots, h_{s_{r-1}+1}\cdot \ldots \cdot h_{s_r}) .     $$
Note that the multiplications are well-defined, since $h_1,\ldots,h_{s_1-1}\in K_{c_1}\subseteq \widetilde{K}_1 = K_{1,s_1}$ and analogously for the other cases.
One then checks that both maps are equivariant with respect to the actions $\chi_1$ and $\chi'_1$.
Therefore we obtain smooth maps on the quotients
$$  \xi  \colon \Gamma_{c_1}\to \Gamma'_{c_1}  \quad \text{and} \quad \zeta\colon \Gamma'_{c_1}\to \Gamma_{c_1} .  $$
It is a direct computation to verify that they are inverses to each other.
\end{proof}
\begin{prop} \label{prop_product_completing_mfld}
The completing manifold $\Gamma_{c}$ is diffeomorphic to the product $\Gamma_{c_1}\times \Gamma_{c_2}$ of the intrinsically defined completing manifolds for $c_1$ and $c_2$, respectively.
Moreover, the diffeomorphism respects the structure as a completing manifold, i.e. the following diagrams commute
$$
\begin{tikzcd}
    \Gamma_c \arrow[]{r}{f}  \arrow[swap]{d}{\cong} & [2.5em] \Lambda M \arrow[]{d}{\cong}
    \\
    \Gamma_{c_1}\times \Gamma_{c_2} \arrow[]{r}{f_1\times f_2} & \Lambda M_1\times \Lambda M_2 
\end{tikzcd}
\hspace{1.5cm}
\begin{tikzcd}
    \Gamma_c \arrow[]{r}{\cong }   \arrow[swap]{d}{p}   & \Gamma_{c_1}\times \Gamma_{c_2} \arrow[]{d}{p_1\times p_2}
    \\
    \Sigma_c \arrow[]{r}{\cong}  & \Sigma_{c_1}\times \Sigma_{c_2}   .
\end{tikzcd}
$$
\end{prop}
\begin{proof}
    This is now a direct consequence of Lemma \ref{lemma_diffeo_gamma_prime} and the definitions of the respective maps.
\end{proof}

\begin{theorem} \label{theorem_intersection_product_completing}
    Let $M_1$ and $M_2$ be compact symmetric spaces and let $\Sigma_1\subseteq\Lambda M_1$ and $\Sigma_2\subseteq \Lambda M_2$ be critical submanifolds consisting of non-trivial closed geodesics.
    We consider the corresponding completing manifolds $\Gamma_1$ and $\Gamma_2$.
    Let $X\in\mathrm{H}_*(\Gamma_1 \times \Gamma_2)$ be a homology class.
    Then the string topology coproduct satisfies
    $$     \vee \, \big((f_1\times f_2)_* X \big) = 0 .      $$
\end{theorem}
\begin{proof}
    We will show that the basepoint intersection multiplicity
    $$    \mathrm{int} \big((f_1 \times f_2)_* X\big) = 1.      $$
    Let $\gamma_1\in\Sigma_1$ and $\gamma_2\in\Sigma_2$.
    There is a time $0<t_a<1$ such that there are no self-intersections of $\gamma_1$ on the interval $(0,t_a]$ and there is a time $0<t_b<1$ such that there are no self-intersections of $\gamma_2$ on the interval $[t_b,1)$.
    It is then clear from the definition of the Ziller cycles that for all $u\in \Gamma_1$ there are no self-intersections in $f_1(u)\in\Lambda M_1$ on the interval $(0,t_a]$.
    Similarly, if $v\in\Gamma_2$ there are no self-intersections in $f_2(v)\in\Lambda M_2$ on the interval $[t_b,1)$.
    We now define maps $\varphi_1\colon \Gamma_1\to \Lambda M_1$ and $\varphi_2\colon \Gamma_2\to \Lambda M_2$.
    Explicitly we set
    $$  \varphi_1(u)(t) = \begin{cases}
        f_1(u)(2t_a \cdot t) &   0\leq t\leq \tfrac{1}{2} \\
        f_1(u)( (2-2t_a) \cdot t + 2 t_a -1) &  \tfrac{1}{2} \leq t \leq 1 
    \end{cases}         $$
     for $u\in\Gamma_1$ and 
    $$  \varphi_2(v)(t) = \begin{cases}
        f_2(v)(2t_b \cdot t) &   0\leq t\leq \tfrac{1}{2} \\
        f_2(v)( (2-2t_b) \cdot t + 2t_b-1 ) &  \tfrac{1}{2} \leq t \leq 1 
    \end{cases} $$
         for $v\in\Gamma_2$.    
         It is clear that $f_1\simeq \varphi_1$ and that $f_2\simeq \varphi_2$ since the maps $\varphi_1$ and $\varphi_2$ are just reparametrizations of the loops appearing in the Ziller cycles.
         If we choose a representing cycle of $X$, i.e. $x\in\mathrm{C}_*(\Gamma_1\times \Gamma_2)$, then we have an inclusion
         $$      \mathrm{im} \big( (\varphi_1\times \varphi_2)_* ( x) \big) \subseteq \mathrm{im}\big( \varphi_1\times \varphi_2 \colon \Gamma_1\times \Gamma_2\to \Lambda M\big) .     $$
        We now show that no loop in the image of $\varphi_1\times \varphi_2$ has a non-trivial self-intersection.
        Assume there was such a loop, i.e. 
        $$   \eta = (\eta_1,\eta_2)  = (\varphi_1(u),\varphi_2(v)) \quad \text{with}\,\,\, u\in\Gamma_1,\,\,v\in\Gamma_2    $$
        and there is a time $t\in (0,1)$ with $\eta(t) = \eta(0)$.
        This implies that 
        $$    \eta_1(t) = \eta_1(0) \quad \text{and} \quad \eta_2(t) = \eta_2(0) .   $$
        But by the definitions of $\varphi_1$ and $\varphi_2$ such a $t$ cannot exist.
        Therefore all loops in the image of $\varphi_1\times \varphi_2$ have no self-intersections.
        This implies that the basepoint intersection multiplicity satisfies $$\mathrm{int}\big( (f_1\times f_2)_* X\big) = 1$$ and therefore the coproduct vanishes by \cite[Theorem 3.10]{hingston:2017}.    
\end{proof}

\begin{cor}
    Let $M_1$ and $M_2$ be compact symmetric spaces and $M  = M_1\times M_2$.
    Take homology with $\mathbb{Z}_2$-coefficients.
    The homology of the free loop space of $M$ relative to the constant loops splits as
    \begin{eqnarray*}
           & & \mathrm{H}_{\bullet}(\Lambda M,M)\cong  \\ &  &\mathrm{H}_{\bullet} (\Lambda M_1,M_1)\otimes \mathrm{H}_{\bullet}(\Lambda M_2,M_2) \,\oplus \mathrm{H}_{\bullet}(\Lambda M_1,M_1)\otimes \mathrm{H}_{\bullet}(M_2) \,\oplus \,   \mathrm{H}_{\bullet}(M_1)\otimes \mathrm{H}_{\bullet}(\Lambda M_2,M_2)  
    \end{eqnarray*}
    and the string topology coproduct is trivial on the first summand $\mathrm{H}_{\bullet} (\Lambda M_1,M_1)\otimes \mathrm{H}_{\bullet}(\Lambda M_2,M_2)$.
    
    The same holds with arbitrary field coefficients if $M_1$ and $M_2$ are orientable and if all critical manifolds as well as all completing manifolds are orientable.
\end{cor}

\bibliography{lit}
 \bibliographystyle{amsalpha}
 
\end{document}